\documentclass[final,12pt]{elsarticle}

\usepackage{graphicx}

\usepackage{amsmath}
\usepackage{amssymb}
\usepackage{textcomp}
\usepackage{float}

\usepackage{color}
\usepackage{units}
\usepackage{geometry}

\journal{Control Engineering Practice}

\newcommand{\nicetilde}{{\raise.17ex\hbox{$\scriptstyle\mathtt{\sim}$}}}
\newcommand{\D}{\mathrm{D}}
\newcommand{\Dha}{\Delta_h^\alpha}
\renewcommand{\d}{\mathrm{d}}
\renewcommand{\Re}{\mathbb{R}}
\newcommand{\N}{\mathbb{N}}
\newcommand{\dfn}{\mathrel{:}=}
\DeclareMathOperator*{\minimize}{minimize}

\begin{document}

\begin{frontmatter}


\title{Model Predictive Control for Offset-Free Reference Tracking of Fractional Order Systems}



\author[NTUA]{Sotiris Ntouskas}
\author[NTUA]{Haralambos Sarimveis}
\author[KUL]{Pantelis Sopasakis}

\address[NTUA]{National Technical University of Athens (NTUA), School of Chemical Engineering, 9~Heroon~Polytechneiou~Street, 15780 Zografou~Campus, Athens, Greece. Fax:~+30-210-7723138, e-mail address:~hsarimv@central.ntua.gr}
\address[KUL]{KU Leuven, Department of Electrical Engineering (ESAT), STADIUS Center for Dynamical Systems, Signal Processing and Data Analytics, Kasteelpark~Arenberg~10, 3001~Leuven, Belgium.}
\begin{abstract}
In this paper an offset-free model predictive control scheme is presented for fractional-order systems using the Gr\"unwald-Letnikov derivative. The infinite-history fractional-order system is approximated by a finite-dimensional state-space system and the modeling error is cast as a bounded disturbance term. Using a state observer, it is shown that the unknown disturbance at steady state can be reconstructed and modeling errors and other persistent disturbances can be attenuated. The effectiveness of the proposed controller-observer ensemble is demonstrated in the optimal administration of an anti-arrhythmic medicine with fractional-order pharmacokinetics.
\end{abstract}

\begin{keyword}
Fractional-order systems \sep Model Predictive Control \sep Gr{\"u}nwald-Letnikov derivative \sep Controlled drug administration \sep Fractional pharmacokinetics.

\end{keyword}
\end{frontmatter}



\section{Introduction}\label{S:1}



\subsection{Background}\label{sec:background}
Fractional calculus is a mathematical extension of the classic calculus of integer-order derivatives and integrals. In fractional calculus, derivatives and integrals are extended to non-integer orders which possess fascinating properties. 
One of the most remarkable properties of fractional-order derivatives is that they are nonlocal operators, that is,  unlike their integer-order counterparts, they cannot be evaluated at a certain point solely by knowing how the function behaves in a neighborhood of this point; instead, the whole history of the function needs to be taken into account~\cite{Podlubny:1999qr}. 

Fractional dynamics have been used to model phenomena exhibiting hereditary properties and long or infinite memory transients. Such phenomena include semi-infinite transmission lines with losses~\cite{Clarke.etal.:2004p}, viscoelastic polymers~\cite{Hilfer:2000qr}, magnetic core coils~\cite{Schafer.etal:2006p}, ultra capacitors~\cite{Gabano201586}, anomalous diffusion in semi-infinite transmission bodies~\cite{Guo.etal:2015p} and several bio-medical applications~\cite{Magin:2010p,Magin.etal:2011p,sopsarmacdok2017jpkpd}. Fractional systems find also several applications in physics~\cite{Hilfer:2000qr}. \citet{Podlubnyetal:2016} offer a thorough overview of the wealth of available toolboxes and software that allow the simulation and controller design for fractional-order systems.
Fractional-order systems and controllers have made their appearance in industrial applications by the extension of the classical PID controller to fractional-order PI$^\lambda$D$^\mu$ ones~\cite{Roy:2016,Beschi2016190,Monje2008798, FeliuBatlle2016159}. In~\cite{FeliuTalegon2016210} fractional order control is applied for the attenuation of vibrations in flexible structures.

During the last few years, a number of works appeared in the literature on the development of Model Predictive Control (MPC) methodologies for fractional order systems.  MPC has gained great popularity in industry and academia due to its inherent capability to take into account state and input constraints, handle complex system dynamics and be resilient to external disturbances~\cite{rawlings2009model}.
In~\cite{Boudjehem:2010,Boudjehem:2012} and~\cite{Romeroetal:2013}, MPC formulations were presented, based on simple input-output fractional order models and using integer-order approximations of the transfer function of the system.  In~\cite{Romeroetal:2013}, the proposed fractional order MPC was demonstrated on the low-speed control of gasoline-propelled cars. In~\cite{Joshietal:2015}, both input-output and state space fractional order models were considered as predictive models in MPC.  In~\cite{Domek:2011} the use of fractional order Takagi-Sugeno fuzzy models was proposed  in the synthesis of fractional MPC. 

A common limitation of all aforementioned works is that although they recognize the importance of MPC in handling input and output/state constraints, they do not take explicitly into account those constraints in the proposed MPC formulations. In~\cite{rhouma2014robust} input constraints were considered in the formulation of the MPC optimization problem, however, without constraints on the state variables and without theoretical stability guarantees.

In practical applications, integer-order derivatives are often used to approximate fractional-order systems~\cite{Vinagre00someapproximations}. Unfortunately these methods come without guarantees of stability and 
satisfaction of constraints. In particular, when these approximations are based on frequency-domain procedures, 
no stability guarantees can be derived and this is a severe shortcoming in safety-critical applications such as 
drug administration. It should be noted again that the behavior of fractional-order systems depends on the whole 
history of their trajectories, therefore it is very difficult to provide cogent evidence based on simulations alone 
without theoretical backup.

An alternative approach has been proposed by~\citet{Guerman:2010qr} where fractional-order
systems are modeled as infinite-dimensional state space systems leading to theoretically interesting stability 
conditions which, nonetheless, are not tractable and cannot be used for controller design.

In our previous work, a controller design approach based on the Gr{\"u}nwald-Letnikov scheme~\cite{Sopasakisetal.:2015jd,SopasakisSarimveis:2017,sopsarmacdok2017jpkpd} was proposed. 
A finite-dimensional approximation was introduced to arrive at a linear time-invariant system and cast 
the discrepancy between the real and the approximate system as an additive bounded uncertainty term. 
A worst-case MPC formulation was presented which leads to \textit{asymptotically stable} behavior 
towards the origin in presence of state and input constraints, even when an approximate finite-history 
model is employed, unlike alternative approaches~\cite{Romero.etal:2009p,Guerman:2010qr}.

\subsection{Contributions}
In this paper an MPC formulation is proposed that achieves offset-free reference tracking for fractional-order systems taking into account input/state constraints. Modeling error is cast as a disturbance term with which the state-space model of the nominal dynamical system is augmented. A state observer is then used to simultaneously estimate the system state and the disturbance using a simple disturbance model. As a result, the closed-loop system can reject disturbances associated with the aforementioned finite-memory approximation, but also other modeling errors due to inexact knowledge of the system parameters, while guaranteeing constraint satisfaction. 

Unlike the controller design approaches we discussed in Section~\ref{sec:background}, which use frequency-domain approximations (i.e., integer-order approximations of the transfer function), the approach we propose in this paper uses a time-domain approximation based on the Gr\"unwald-Letnikov derivative.

The proposed MPC strategy is demonstrated in a case study emerging from pharmacokinetics and pharmacodynamics, where fractional-order systems are becoming increasingly popular over the last years. The work of~\citet{Kytariolos.etal:2010p} introduced fractional-order dynamics in pharmacokinetics, highlighting why the classical \textit{in-vitro-in-vivo correlations} theory fails. Certain nonlinearities, anomalous diffusion, deep tissue trapping, diffusion across fractal manifolds such as systems of capillaries, synergistic and competitive actions are cases that can hardly be modeled by integer-order systems~\cite{DokMach:2008p}. Fractional-order pharmacokinetic dynamics can be cast as physiologically based pharmacokinetic (PBPK) or compartmental models by properly re-writing the mass balance equations using fractional-order derivatives in a way that mass balances are not violated~\cite{Dokoumetzidis:2010a}. The controlled drug administration for drugs with fractional dynamics is a key enabler of an effective and realistic therapy and a valuable tool for the clinical practice, especially in presence of constraints~\cite{SopasakisSarimveis:2014}.

This paper is organized as follows:  In Section~\ref{S:2} we describe fractional-order dynamical systems in terms of the Gr{\"u}nwald-Letnikov derivative and derive control-oriented approximations of bounded error. In Section~\ref{S:3} we propose an MPC scheme for offset-free control using a state observer for an augmented system which is able to attenuate modeling errors. Lastly, In Section~\ref{S:4} we present such an offset-free MPC for the control of amiodarone administration to patients and show that the controlled system is resilient to inexact knowledge of the pharmacokinetic parameters of the patients (which are, typically, not known).

\subsection{Notation}
\label{sec:notation}

Hereafter, $\mathbb{R}$ and $\mathbb{N}$ denote the sets of real and nonnegative 
integers respectively. We denote by $\mathbb{N}_{[k_1, k_2]}$ the set of all integers 
in the closed interval $[k_1, k_2]$. The set of real $n$-dimensional vectors is 
denoted by $\Re^n$ and the set of $m$-by-$n$ matrices by $\Re^{m\times n}$.
All sets are denoted by calligraphic uppercase letters and all matrices
are denoted by uppercase letters. Vectors and scalars are denoted by lowercase letters. 
The transpose of a matrix $A$ is denoted by $A'$.

\section{Fractional-Order Systems}\label{S:2}
\subsection{Discrete-time fractional operators}
In this section a fractional-order differential operator, the Gr{\"u}nwald-Letnikov derivative is introduced. Let $f:\Re\to\Re^n$ be a bounded function. Let us first introduce the \textit{Gr{\"u}nwald-Letnikov difference} of $f$ at $t$ of order $\alpha > 0$ and step size $h>0$, which is defined as
\begin{align}
\label{eq:GrunwLetBack}
\Dha f(t) = \sum_{j=0}^{\infty}(-1)^j \tbinom{\alpha}{j} f(t-j h)
\end{align}
where $\binom{\alpha}{0}=1$ and for $j \in \mathbb{N}, j>0$, $\tbinom{\alpha}{j} = \prod_{i=0}^{j-1}\frac{\alpha-i}{i+1}$.
Furthermore, let us define $c_j^\alpha = (-1)^j\tbinom{\alpha}{j}$ and note that $|c_j^\alpha| \leq \alpha^j/j!$ 
for all $j \in \mathbb{N}$, therefore, the sequence $c_j^\alpha$ is absolutely summable and $\Dha$ is well defined. 
It is now clear that in order to estimate $\Dha f(t)$ for non-integer orders $\alpha$ the whole history of $f(t)$ is needed.  

The Gr{\"u}nwald-Letnikov operator leads to the definition of the Gr{\"u}nwald-Letnikov fractional derivative of order $\alpha$ as
\begin{align}
\label{eq:GLderivative}
\D^{\alpha}f(t) = \lim_{h\rightarrow 0}\frac{\Dha f(t)}{h^{\alpha}} 
\end{align}
provided that the limit exists. It can be verified that $\D^k$, for $k\in\N$, are the ordinary integer-order derivatives with respect to $t$ and, by convention, $\D^0 f(t) = f(t)$.

Using the definition above, it is easy to describe fractional-order dynamical systems with state $x: \mathbb{R} \rightarrow \mathbb{R}^n$ and input $u: \mathbb{R} \rightarrow \mathbb{R}^m$ as 
\begin{align}
\label{eq:FracSysDescription}
\sum_{i=0}^l A_i \D^{\alpha_i}x(t) = \sum_{i=0}^r B_i \D^{\beta_i}u(t)
\end{align}
where $l, r \in \mathbb{N}$, $A_i$ and $B_i$ are matrices of appropriate dimensions and all powers $\alpha _i$ and $\beta_i$ are nonnegative.

For the discretization of a fractional system, an Euler-type method is used to approximate $\D^{\alpha}$ in~\eqref{eq:FracSysDescription}, for a fixed time step size $h$, using $h^{-\alpha} \Dha$. Using the forward operator for the derivatives of the states, and the backward operator for the input variables, the discretization of Equation~\eqref{eq:FracSysDescription} becomes 
\begin{align}
\label{eq:FracSysDDescription}
\sum_{i=0}^l \bar{A}_i\, \Delta_h^{\alpha_i}x_{k+1} = \sum_{i=0}^r \bar{B}_i \,\Delta_h^{\beta_i}u_k
\end{align}
For convenience in~\eqref{eq:FracSysDDescription} it is: $x_k = x(kh), u_k = u(kh)$ for $k \in \mathbb{Z}$, and $\bar{A}_i = h^{-\alpha _i} A_i, \bar{B}_i = h^{-\beta _i} B_i$.

As it becomes obvious from the equations above, in order to fully describe a fractional-order system, infinite-dimensional operators of the form $\Delta_h^\alpha$ should be used. Therefore, it is practically impossible to simulate such systems or design feedback controllers. In what follows, Equation~\eqref{eq:FracSysDDescription} will be approximated using the methodology described in~\cite{Sopasakisetal.:2015jd} where a finite-dimensional approximation is proposed and the approximation error is treated as a bounded additive disturbance.

\subsection{Bounded-error Finite-dimension Approximation}\label{sec:fin-dim-approx}
As infinite-memory systems pose severe limitations regarding their simulation and controller design, the \textit{truncated Gr{\"u}nwald-Letnikov difference operator} of length $\nu$ is introduced
\begin{align}
\label{eq:truncatedGLdiff}
\Delta^{\alpha}_{h,\nu}x_k = \sum_{j=0}^{\nu}c_j^{\alpha}x_{k-j},
\end{align}
and the approximation of~\eqref{eq:FracSysDescription} for $\nu \geq 1$ becomes
\begin{align}
\label{eq:approxSystem}
\sum_{i=0}^{l}\bar{A}_i\Delta^{\alpha_i}_{h,\nu}x_{k+1} = \sum_{i=0}^{r}\bar{B}_i\Delta^{\beta_i}_{h,\nu-1}u_{k},
\end{align}
which is equivalently written as 
\begin{align}
\label{eq:approxSystemSS}
\sum_{j=0}^{\nu}\hat{A}_j x_{k-j+1} = \sum_{j=0}^{\nu}\hat{B}_j u_{k-j}
\end{align}
where $\hat{A}_j = \sum_{i=0}^{l} \bar{A}_i c_j^{\alpha_i}$ and $\hat{B}_j = \sum_{i=0}^{r} \bar{B}_i c_j^{\beta_i}$ for $j \in \mathbb{N}_{[0,\nu]}$. Assuming that matrix $\hat{A}_0$ is nonsingular,  $\tilde{A}_j = -\hat{A}_0^{-1} \hat{A}_{j+1}$ and $\tilde{B}_j = \hat{A}_0^{-1} \hat{B}_j$ are defined and the last Equation~\eqref{eq:approxSystemSS} becomes
\begin{align}
\label{eq:approxSystemSS2}
x_{k+1} = \sum_{j=0}^{\nu-1}\tilde{A}_j x_{k-j} + \sum_{j=1}^{\nu}\tilde{B}_j u_{k-j} + \tilde{B}_0 u_k
\end{align}
Considering as the state variable $\tilde{x}_{k}=(x_k, x_{k-1},...,x_{k-\nu+1}, u_{k-1},...,u_{k-\nu})'$, the system can be written in the compact form of a finite-dimensional linear time-invariant (LTI) system as
\begin{align}
\label{eq:approxSystemSS3}
\tilde{x}_{k+1} = A\tilde{x}_k + Bu_k.
\end{align}
In the next section, this form will be used to formulate a model predictive control problem, assuming that the pair $(A,B)$ is stabilizable. 

The fractional-order difference operator can be written as the sum of the truncated difference operator and a \textit{residual} term containing all those past states that were not taken into account, that is
\begin{align}
\label{eq:GLoptruncAndRes}
\Dha = \Delta_{h,\nu}^{\alpha} + R_{h,\nu}^{\alpha},
\end{align}
where $R_{h,\nu}^{\alpha}$ is the operator 
\begin{align}
\label{eq:Residuals}
R_{h,\nu}^{\alpha}(x_k) = h^{-\alpha} \sum_{j{}={}\nu{}+{}1}^{\infty}\!c_{j}^{\alpha}x_{k-j}.
\end{align}
Assume that up to time $k$ the state has been constrained in a compact convex set $\mathcal{X}\subseteq \Re^n$ containing the origin in its interior, that is $x_{k-j} \in \mathcal{X}$ for all $j \in \mathbb{N}$. Let us denote the \textit{Minkowski sum} of two sets $\mathcal{C},\mathcal{D}\subseteq \Re^n$ as the set $\mathcal{C}\oplus \mathcal{D} \dfn \{c+d: c\in \mathcal{C}, d\in \mathcal{D}\} \subseteq \Re^n$. For a collection of sets $\mathcal{C}_1,\mathcal{C}_2,\ldots$, their Minkowski sum $\mathcal{C}_1\oplus \mathcal{C}_2 \oplus \cdots$ is denoted as $\bigoplus_{k=1}^{\infty} \mathcal{C}_k$.

By virtue of~\eqref{eq:Residuals},  
\begin{align}
\label{eq:ResidualsSet}
R_{h,\nu}^{\alpha}(x_k) \in \bigoplus_{j{}={}\nu{}+{}1}^{\infty}\!c_{j}^{\alpha} \mathcal{X}.
\end{align}
In the special case where $\mathcal{X}$ is also a \textit{balanced} set, that is $-x\in \mathcal{X}$ for all $x\in\mathcal{X}$, \eqref{eq:ResidualsSet} reduces to
\begin{align}
\label{eq:ResidualsSet-balanced}
R_{h,\nu}^{\alpha}(x_k) \in \left({}\sum_{j{}={}\nu{}+{}1}^{\infty}|c_j^{\alpha}|{}\right)\!\mathcal{X}.
\end{align}
In this case, it is evident that $\nu$ directly controls the size of the right hand side of~\eqref{eq:ResidualsSet-balanced}. In particular, $({}\sum_{j=\nu_1+1}^{\infty}|c_j^{\alpha}|{})\mathcal{X} \supset ({}\sum_{j=\nu_2+1}^{\infty}|c_j^{\alpha}|{})\mathcal{X} $ for all $\nu_1 < \nu_2$, that is, the larger the value of $\nu$ is, the lower the worst-case approximation error will be. 

The approximation error can be integrated into the state space form as an additive disturbance term, so~\eqref{eq:approxSystemSS2} becomes
\begin{align}
\label{eq:UncertainSystem}
\tilde{x}_{k+1}{}={}A\tilde{x}_k{}+{}Bu_k{}+{}Gd_k
\end{align}
with $G = \begin{bmatrix}
I & 0 & \cdots & 0 \\
\end{bmatrix}'$. 
The disturbance term $d_k$ is introduced to cast the approximation error which is in fact a bounded error term as we are about to show.
Assuming that the state is constrained, as discussed above, in a set $\mathcal{X}$ and, similarly, the input is constrained in a convex compact set $\mathcal{U}$ containing the origin in its interior, following~\cite{Sopasakisetal.:2015jd}, $d_k$ is bounded in a compact set $\mathcal{D}_{\nu}$ given by
\begin{align}
\label{eq:DistSet}
\mathcal{D}_{\nu}{}={}\mathcal{D}_{\nu}^x{}\oplus{}\mathcal{D}_{\nu}^u
\end{align}
where sets $\mathcal{D}_{\nu}^x$ and $\mathcal{D}_{\nu}^u$ are\footnote{Simpler formulas for $\mathcal{D}_{\nu}^x$ and $\mathcal{D}_{\nu}^u$  may be derived assuming that $\mathcal{X}$ and $\mathcal{U}$ are balanced sets, i.e., they are symmetric about the origin --- or equivalently $-x\in \mathcal{X}$ whenever $x\in \mathcal{X}$ and $-u\in \mathcal{U}$ whenever $u\in \mathcal{U}$. The reader is referred  to~\cite{Sopasakisetal.:2015jd} for further information.}
\begin{subequations}
\label{eq:DistSetTypes}
\begin{align}
\mathcal{D}_{\nu}^x &= \bigoplus_{i{}={}0}^l -\hat{A}_0^{-1} \bar{A}_i \bigoplus_{j{}={}\nu{}+{}1}^{\infty}c_j^{\alpha_i} \mathcal{X},\\
\mathcal{D}_{\nu}^u &= \bigoplus_{i{}={}0}^l \hat{A}_0^{-1} \bar{B}_i \bigoplus_{j{}={}\nu{}+{}1}^{\infty}c_j^{\beta_i} \mathcal{U}.
\end{align}
\end{subequations}
In the special case where $\mathcal{X}$ and $\mathcal{U}$ are balanced sets, these expressions can be further simplified as explained above.

The assumption that $\mathcal{X}$ and $\mathcal{U}$ contain the origin in their interiors is standard in regulation problems, where the objective is to steer the state of the system to the origin.  In case of reference tracking, it is desirable that the system output $y_k$, which is typically a linear combination of states $y_k=Cx_k$, converges to a certain set-point $r$. Assume that the set-point is attainable by a feasible state $\bar{x}\in \mathcal{X}$ so that $C\bar{x}=r$. It then suffices that $x_k$ converges to $\bar{x}$ with $x_k\in \mathcal{X}$ with a feasible sequence of control actions $u_k\in \mathcal{U}$ which converges to a feasible point $\bar{u}\in \mathcal{U}$. Then, the deviation variables $\delta x_k = x_k -\bar{x}$ and $\delta u_k = u_k - \bar{u}$ may be defined. The above convergence requirements become $\delta x_k \to 0$ and $\delta u_k \to 0$ with constraints $\delta x_k\in \{-\bar{x}\}\oplus \mathcal{X}$ and $\delta u_k \in \{-\bar{u}\}\oplus \mathcal{U}$. In order to be able to derive bounded-error approximations of the fractional order system with input $\delta u_k$ and state $\delta x_k$ it needs to be assumed that $\bar{x}$ lies in the interior of $\mathcal{X}$ and $\bar{u}$ lies in the interior of $\mathcal{U}$.

Larger memory lengths lead to better approximations since $R_{h,\nu}^{\alpha}(x)\to 0$
as $\nu\to\infty$, that is, arbitrarily low approximation errors can be obtained by choosing 
an appropriately large value of $\nu$. 
That said, on one hand it is desirable to opt for a large value of $\nu$ 
to obtain better approximations. On the other hand, too large values of $\nu$ should be 
avoided as they will incur a larger computation cost especially regarding the solution of 
the MPC problem in real time (see Section~\ref{sec:offset-free-mpc}).

\section{Offset-free reference tracking}\label{S:3}

\subsection{Observer design}\label{sec:obsv-design}
In this section, it is assumed that the output of~\eqref{eq:UncertainSystem} can be measured, which is given by
\begin{align}
 y_k = C\tilde{x}_k + C_d d_k.
\end{align}
A state observer will be designed to simultaneously reconstruct the state $\tilde{x}_k$ and persistent disturbances $d_k$ using the simple disturbance model $d_{k+1} = d_k$. 
In~\cite{Rajamanietal:2009} it is shown that even if the disturbance model is not accurate, the offset-free
regulation properties of the closed-loop system and the closed-loop performance will not be altered.

We introduce an \textit{augmented state observer} which is a dynamical system which using input-output measurements, will produce estimates of the augmented state variable $\xi_k = \big[\tilde{x}_k', d_k'\big]'$ whose dynamics is
\begin{subequations}\label{eq:Observer1}
\begin{align}
\xi_{k+1} &= \bar{A} \xi_k + \bar{B} u_k, \\
y_k &= \bar{C} \xi_k,
\end{align}
\end{subequations}
with
\begin{align}
\label{eq:Observer2}
\bar{A} = 
  \begin{bmatrix}
    A & G \\ 
    0 & I
  \end{bmatrix}, 
\bar{B} = 
  \begin{bmatrix}
    B\\ 
    0
  \end{bmatrix}, 
\bar{C} = 
  \begin{bmatrix}
    C & C_d 
  \end{bmatrix}
\end{align}
where $I$ is the identity matrix. 
According to~\cite[Prop.~1]{Maeder:2009p}, system~\eqref{eq:Observer1} is observable if and only if $(C,A)$ is observable and 
\begin{align*}
 \begin{bmatrix}
  A-I & G\\
  C   & C_d
 \end{bmatrix},
\end{align*}
is full column rank. For these conditions to be satisfied, we need to choose the dimension of $d_k$
to be no larger than the number of outputs. We shall hereafter assume that $d_k$ has the same dimension as $y_k$.
Additionally, matrix $C_d$ reflects the effect that the persistent disturbance $d_k$ has (directly) on the system output.

The observer is a dynamical system that produces estimates $\hat{\xi}_k = \big[\hat{\tilde{x}}_k', \hat{d}_k'\big]'$,  that is, the observer produces simultaneously estimates of the system state and the disturbance. Here we use a linear observer of the form
\begin{align}
\label{eq:Luenberger}
\hat{\xi}_{k+1} &= \bar{A}\hat{\xi}_k + \bar{B}u_k + Le_k,\\
e_k &= \bar{C}\hat{\xi}_k - y_k,
\end{align}
where $e_k$ is the state estimation error. 

By properly choosing the gain matrix $L$, either by pole placement methods or with LQG methods, it can be guaranteed that the estimation error converges to $0$ for any initial estimate $\hat{\xi}_0$. The estimates of the augmented state are then provided to the tracking MPC which is described in the following section.

\subsection{Offset-free model predictive control}\label{sec:offset-free-mpc}
In MPC, the control action is computed at every time instant by minimizing an index which quantifies the performance of the system along a finite prediction horizon. This performance index is used to choose an optimal sequence of control actions among the set of such admissible sequences, while corresponding state sequences are predicted using a system model. The first element of the optimal sequence is then applied to the system.

\citet{Maeder:2009p} proposed an offset-free MPC problem formulation without integral action where the MPC problem makes use of the augmented model~\eqref{eq:Observer1}. The main concept of this approach is that the observer eventually reconstructs any persistent disturbances and this information is used to steer the state to the desired set-point ${r_k}$ without offset. The objective of the closed-loop system is for the tracking error $y_k-r_k$ to converge to $0$ whenever $r_k$ is a convergent sequence.

At every time instant $k$ the following finite-horizon optimal control problem is solved:
\begin{align}
\label{eq:ofMPC}
\mathsf{P}_{\mathrm{mpc}}(\hat{\tilde{x}}_k, \hat{d}_k): \minimize_{\textbf{u}=\{u_{k+j\mid k}\}_{j=0}^{N-1}} V_N(\textbf{u}; \bar{x}_k, \bar{u}_k)
\end{align}
subject to the constraints
\begin{subequations}
\label{eq:Constraints}
\begin{align}
\tilde{x}_{k+j+1\mid k} &= A\tilde{x}_{k+j\mid k} + Bu_{k+j\mid k} + G d_{k+j\mid k}, j \in \mathbb{N}_{[0, N-1]} \label{eq:EqConstraints-1}\\
d_{k+j+1\mid k} &= d_{k+j\mid k}, j \in \mathbb{N}_{[0, N-1]}\label{eq:EqConstraints-2}\\
F_x \tilde{x}_{k+j\mid k} &+ F_u u_{k+j\mid k} \leq f, j \in \mathbb{N}_{[0, N-1]}\label{eq:EqConstraints-3}\\
\tilde{x}_{k\mid k} &= \hat{\tilde{x}}_k \label{eq:EqConstraints-4}\\
d_{k\mid k} &= \hat{d}_k\label{eq:EqConstraints-5}
\end{align}
\end{subequations}
with cost function
\begin{align}\label{eq:CostFunction}
V_N( \mathbf{u}; \bar{x}_k, \bar{u}_k) &\dfn
 \|\tilde{x}_{k+N\mid k}-\bar{x}_{k}\|_P^2 \notag\\
   &+ \sum_{j=0}^{N-1} 
         \left(
            \|\tilde{x}_{k+j\mid k}-\bar{x}_k\|_Q^2 + \|u_{k+j\mid k}-\bar{u}_k\|_R^2
         \right) 
\end{align}	
where the notation $\|x\|_F^2 = x'Fx$ is used. 
The optimization problem~\eqref{eq:ofMPC} is solved over sequences of future control actions $\textbf{u}=\{u_{k+j\mid k}\}_{j=0}^{N-1}$ of length $N$ and the problem is solved using information that is available at time $k$ which is provided by the augmented state observer. The first element of this sequence, $u_{k\mid k}^\star$, is then applied to the system in a \textit{receding horizon} fashion~\cite{rawlings2009model}.

The optimization problem $\mathsf{P}_{\mathrm{mpc}}(\hat{x}_k, \hat{d}_k)$ is solved with respect to the dynamics of the augmented system described by~\eqref{eq:EqConstraints-1} and~\eqref{eq:EqConstraints-2}, state and input constraints~\eqref{eq:EqConstraints-3} and the initial conditions~\eqref{eq:EqConstraints-4} and~\eqref{eq:EqConstraints-5}. The constraints~\eqref{eq:EqConstraints-3}, where $F_x$ and $F_u$ are matrices of appropriate dimensions and $\leq$ is meant element-wise, encompasses simple bounds of the form $\tilde{x}_{\min} \leq \tilde{x}_k \leq \tilde{x}_{\max}$ and $u_{\min} \leq u_k \leq u_{\max}$ or more involved polyhedral joint state-input constraints.

The terms $\bar{x}_k$ and $\bar{u}_k$ in~\eqref{eq:ofMPC} are calculated at runtime by solving the linear system
\begin{align}\label{eq:xhatuhat}
\begin{bmatrix}
  A-I & B \\ 
  C & 0
\end{bmatrix} 
\begin{bmatrix}
  \bar{x}_k \\ 
  \bar{u}_k
\end{bmatrix} = 
\begin{bmatrix} 
  -G \hat{d}_k \\ 
  r_k - C_d \hat{d}_k
\end{bmatrix}.
\end{align}
Equation~\eqref{eq:xhatuhat} is used at every time instant to determine $\bar{x}_k$ and $\bar{u}_k$  given the desired set-point $r_k$ and the disturbance estimate $\hat{d}_k$. In fact, it suffices to find a matrix $W$ so that 
\begin{align}
 \begin{bmatrix}
  A-I & B \\ 
  C & 0
\end{bmatrix} W = 
\begin{bmatrix}
  -G & 0\\
  -C_d & I
\end{bmatrix}.
\end{align}
Then, $\bar{x}_k$ and $\bar{u}_k$ can be obtained by 
\begin{align}
 \begin{bmatrix}
  \bar{x}_k \\ 
  \bar{u}_k
\end{bmatrix} = 
W
\begin{bmatrix}
 \hat{d}_k\\ r_k
\end{bmatrix}.
\end{align}

As the linear observer produces a sequence of disturbances $\hat{d}_k$ which converges to a steady state value $\hat{d}_\infty=d_{\infty}$, conditions~\eqref{eq:xhatuhat} enforce that the tuple ($\bar{x}_\infty, d_{\infty}, \bar{u}_{\infty},r_{\infty})$ is an equilibrium point of~\eqref{eq:Observer1}, i.e., as $r_k\to r_\infty$, $x_k\to\bar{x}_\infty$, $u_k\to \bar{u}_\infty$, we have $y_k\to r_\infty$.

Finally, the positive semidefinite matrix $Q$ and the positive definite matrix $R$ are the tuning knobs of MPC and are typically chosen to be diagonal matrices; larger values of $Q$ lead to a faster convergence of the system states to their equilibria, whereas larger values of $R$ lead to a smoother actuation and slower convergence towards the steady state.
Matrix $P$, which defines the terminal cost function in~\eqref{eq:CostFunction}, is taken to be the unique solution of the Riccati equation
\begin{align}
\label{eq:Riccati}
P = A'PA - A'PB(B'PB+R)^{-1}B'PA + Q.
\end{align}
By choosing $P$ to be the solution of~\eqref{eq:Riccati}, 
the value function of~\eqref{eq:ofMPC} becomes a Lyapunov function for the 
closed-loop system and stability is guaranteed~\cite{Maeder:2009p,rawlings2009model}.

The optimization problem $\mathsf{P}_{\mathrm{mpc}}$ is a convex quadratic optimization 
problem with polyhedral constraints which can be solved very efficiently in practice.
In fact, there exist algorithms which exhibit very fast convergence and scale linearly with the prediction 
horizon~\cite{patrinos2012accelerated,GPAD_NMPC}.

\section{Application: Controlled drug administration}\label{S:4}

\subsection{Problem statement: objectives and constraints}
Fractional-order dynamics are of great interest in applications of pharmacokinetics such as controlled drug administration~\cite{Dokoumetzidis:2010a,Dok.etal:2011p,SopasakisSarimveis:2014}. In this section the reference tracking methodology presented in Section~\ref{S:3} is applied for the controlled intravenous administration of a drug which exhibits fractional-order dynamics.

\begin{figure}[H]
\centering
\includegraphics[width=0.8\columnwidth]{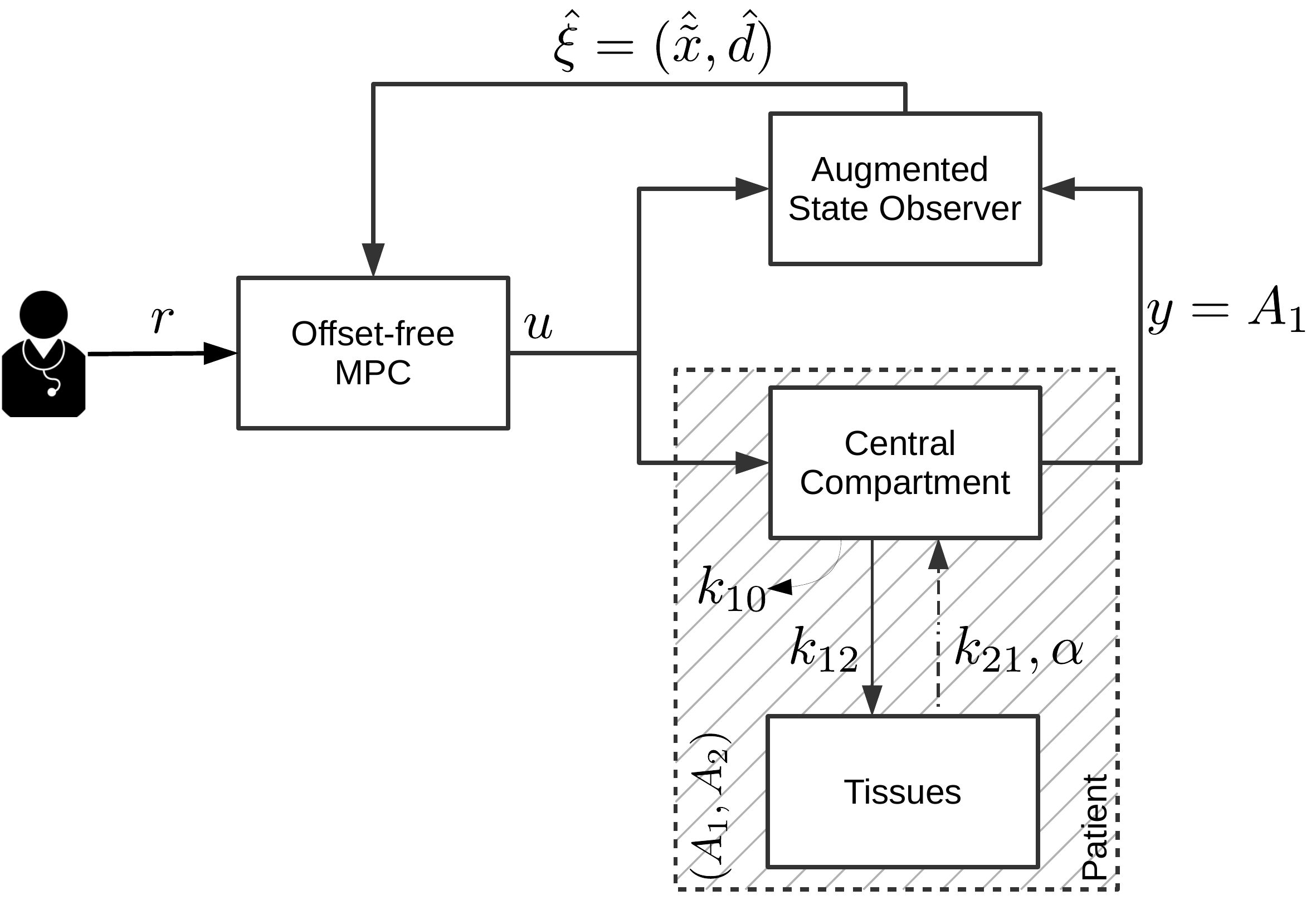}
  \caption{Schematic of controlled drug administration. 
  The pharmacokinetic system is described by a two-compartment model
  from which measurements are available from the central compartment only.
  An augmented state observer produces estimates of the system state $\hat{\tilde{x}}$ 
  and the disturbance $\hat{d}$ which are fed to the controller. The treating 
  physician specifies the desired reference signal $r$.}
    \label{fig:system}
\end{figure}

Amiodarone (CAS registration number: 1951-25-3) is an anti-arrhythmic medication used for the treatment of ventricular tachycardia, shock-refractory ventricular and atrial fibrillation, which has been successfully modeled by a fractional compartmental system in~\cite{Dokoumetzidis:2010a}. The pharmacokinetic model considers two compartments: the central compartment which corresponds to the blood stream where the drug is introduced with rate $u$ and a peripheral compartment which corresponds to the tissues where the blood is distributed. The distribution from the blood stream to the tissues is assumed to follow a simple first-order dynamics with constant $k_{12}$. Amiodarone is excreted from the central compartment following first-order excretion kinetics with rate constant $k_{10}$. A mass flow from the tissues back to the central compartment takes place via anomalous diffusion and follows fractional-order kinetics of order $1-\alpha$ with rate $k_{21}$. This structure is shown in Figure~\ref{fig:system} along with the observer and the offset-free controller.

The mass balance equations, taking into account the fractional-order kinetics, lead to the dynamical system
\begin{subequations}
\label{eq:DokFracSS}
\begin{align}
\frac{\d A_1}{\d t} &= -(k_{12} + k_{10})A_1 +  k_{21}\cdot \D^{1-\alpha}A_2 + u\\
\frac{\d A_2}{\d t} &= k_{12}A_1 - k_{21}\cdot \D^{1-\alpha}A_2,
\end{align}
\end{subequations}
where $A_1$ and $A_2$ are the amounts of amiodarone (in $\unit{ng}$) in the central and peripheral compartments respectively and $u$ is the intravenous administration rate into the central compartment (in $\unitfrac{ng}{day}$). The (nominal) values of the parameters of the system are $\alpha = 0.587$, $k_{10}=\unit[1.4913]{day^{-1}}$, $k_{12} = \unit[2.9522]{day^{-1}}$ and $k_{21} = \unit[0.4854]{ day^{-\alpha}}$. It is considered that amiodarone is administered to the patient intravenously and continuously. Only measurements of $A_1$ are available in real time.

As shown in Figure~\ref{fig:system}, a state observer receives measurements of the system output, that is $y=A_1$ and input variables and produces estimates of the full state $\hat{\tilde{x}}$ and disturbance $\hat{d}$ which are provided to the MPC. The treating physician prescribes the desired set-point $r$ for $A_1$ to the MPC controller which decides the administration rate $u$ at every time instant. 

The continuous-time fractional-order system is simulated using the \textit{Oustaloup filter} --- an integer-order approximation of the fractional-order differentiators $s^a$ in a frequency range $[\omega_L, \omega_H]$ via a logarithmically spaced sampling of $N_f$ points~\cite{OusLev+00,monje2010fractional}. Using the Oustaloup filter, $s^a$ is approximated by a transfer function with numerator and denominator degrees equal to $N_f$.

Oustaloup's filter is chosen because, as shown in~\cite{ifac-fractional-2017}, it is very accurate 
for this dynamical system. Following~\cite{ifac-fractional-2017}, the values $\omega_L=\unit[10^{-2}]{day^{-1}}$, 
$\omega_H=\unit[10^3]{day^{-1}}$ and $N_f=8$ are chosen which are known to lead to a highly accurate 
approximation. Other frequency-domain approximation methods are available in the literature such as the 
Pad\'e and the Matsuda-Fujii which are studied in~\cite{ifac-fractional-2017} in terms of accuracy and 
seem to lead to solutions of comparable quality. We should underline, however, that frequency domain methods
do not come with bounded error guarantees, that is, we cannot know the maximum approximation error in advance.
This motivates the use of the \textit{Gr{\"u}nwald-Letnikov} approximation which is the most suitable approach.
Frequency-domain approximations lead to approximation errors which --- although in practice are small --- 
are not accompanied by known bounds. 
On the other hand, the truncated Gr{\"u}nwald-Letnikov approximation, as less parsimonious as it might be, 
leads to a bounded-error approximation in the time domain which is necessary for the design of MPC controllers 
with theoretical stability and constraint satisfaction guarantees.
Other methods, such as the numerical inverse Laplace transform, lead to highly accurate solutions, but these 
are not suitable for controller design.

In practice people have used other approximations to design MPC controllers 
(than the Gr{\"u}nwald-Letnikov discrete-time derivative)~\cite{Joshietal:2015,RhoBouBou14,rhouma2014robust}, 
however, without any guarantees of constraint satisfaction which is particularly important 
in fail-sensitive applications (e.g., medical).

\subsection{Control-oriented modeling}
The fractional-order pharmacokinetic system~\eqref{eq:DokFracSS} can be written as 
\begin{align}
\label{eq:sysSS}
\D\begin{bmatrix}A_1 \\ A_2\end{bmatrix} = \begin{bmatrix}-k_0 & 0 \\ k_{12} & 0\end{bmatrix} \begin{bmatrix}A_1 \\ A_2\end{bmatrix} + \begin{bmatrix}0 & k_{21} \\ 0 & -k_{21}\end{bmatrix} \D^{\beta} \begin{bmatrix}A_1 \\ A_2\end{bmatrix} + \begin{bmatrix}1 \\ 0\end{bmatrix}  u
\end{align}
with $\beta = 1 - \alpha $ and $k_0 = k_{12}+k_{10}$. 
The constants
\begin{align}
\label{eq:VarDecl}
M = \begin{bmatrix}-k_0 & 0 \\ k_{12} & 0\end{bmatrix}, \ \Theta = \begin{bmatrix}0 & k_{21} \\ 0 & -k_{21}\end{bmatrix}, \ B = \begin{bmatrix}1 \\ 0\end{bmatrix},
\end{align}
are defined and the variable
\begin{align}
 x = \begin{bmatrix}A_1 \\ A_2\end{bmatrix}.
\end{align} is introduced.
Now~\eqref{eq:sysSS} can be rewritten in a compact form as
\begin{align}
\label{eq:compactSS}
\D x = M x + \Theta \D^{\beta}x + Bu.
\end{align}
The above system can now be discretized with sampling time $h$ as discussed in Section~\ref{sec:fin-dim-approx} to yield
\begin{align}
\label{compactSS2}
\frac{x_{k+1}-x_k}{h} = Mx_k + \Theta \Delta_{h,\nu}^{\beta} x_k + B u_k + B_d d_k,
\end{align}
where $B_d = [\begin{smallmatrix}1\\0\end{smallmatrix}]$ and $d_k\in\Re$ is a generic disturbance term which is used to encompass the effect that the following factors have on $A_1$: (i) the approximation error due to the use of the truncated operator $\Delta_{h,\nu}^{\beta}$,  (ii) modeling errors due to inexact knowledge of the actual pharmacokinetic parameters, (iii) modeling errors due to pharmacokinetic dynamics not captured by the model described in the the previous section. We shall now use Equation~\eqref{compactSS2} to derive a simplified linear time-invariant model with a scalar disturbance term which will be used to design the MPC controller and the observer as discussed above.

Let us define the state variable $\tilde{x}_k = (x_k, x_{k-1},\ldots, x_{k-\nu})\in\Re^{2(\nu+1)}$. 
By defining $\Lambda = Mh + I + \Theta h^{\beta}$ the approximate system in a matrix form is derived:
\begin{subequations}\label{eq:ApproxSS:sys}
\begin{align}
\label{eq:ApproxSS}
\tilde{x}_{k+1}
= 
\begin{bmatrix} 
  \Lambda & \Theta h^{\alpha}c_1^{\beta} & \cdots & \Theta h^{\alpha}c_{\nu-1}^{\beta} & \Theta h^{\alpha} c_{\nu}^{\beta} \\
I & 0 & \cdots & 0 & 0 \\
0 & I & \cdots & 0 & 0 \\
\vdots & \vdots & \ddots & 0 & 0 \\
0 & 0 & \cdots & I & 0 
\end{bmatrix}
\tilde{x}_{k} 
+ 
\begin{bmatrix}
Bh\\ 
0\\ 
0\\ 
\vdots\\ 
0 
\end{bmatrix} u_k
+ G d_k,
\end{align}
with $G = \begin{bmatrix}B_d' & 0 & \ldots & 0\end{bmatrix}'$.
The system output is
\begin{align}
 y_k = \begin{bmatrix}1&0&\ldots&0\end{bmatrix}\tilde{x}_k,
\end{align}
\end{subequations}
therefore, $C=\begin{bmatrix}1&0&\ldots&0\end{bmatrix}$ and 
$C_d = 0$. One may also verify that the observability assumptions stated
in Section~\ref{sec:obsv-design} are satisfied for system~\eqref{eq:ApproxSS:sys}.

In what follows the initial condition is $A_1(t)=A_2(t)=0$ for all $t\leq 0$ is assumed, that is $\tilde{x}_0=0$.

\subsection{Controlled drug administration}
Hereafter a discretization with $h = \unit[0.1]{day}$ and a memory 
length $\nu=25$ which corresponds to $2.5$ days are used.
In the closed-loop simulation, we assume that the treating 
physician sets the reference to the system output ($A_1$) equal to  $\unit[0.5]{ng}$ for the first $80$ days 
and increases it to $\unit[1.0]{ng}$ for the next $70$ days.
The administration rate is constrained in the interval 
$0 \leq u_k \leq u_{\max} = \unitfrac[2]{ng}{day}$, while the system output should not exceed the upper bound of $\unit[1.03]{ng}$  . 
The prediction horizon for the controller is set to $N = 60$. 
The weight matrices of the MPC controller are chosen to be $Q=0.25\cdot I$ and $R=5$.

\begin{figure}[ht]
\centering
\includegraphics[width=0.325\columnwidth]{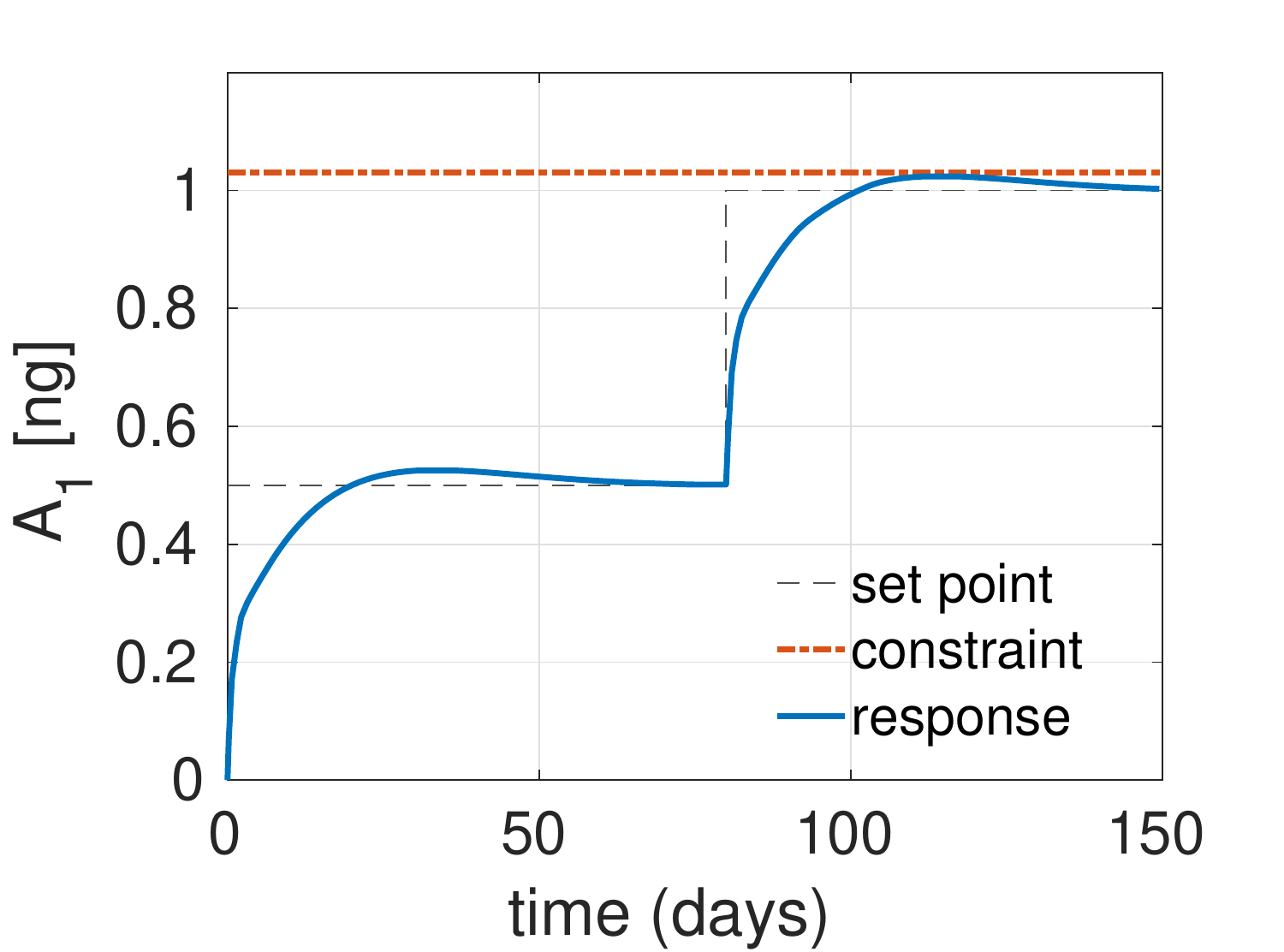}{}%
\includegraphics[width=0.325\columnwidth]{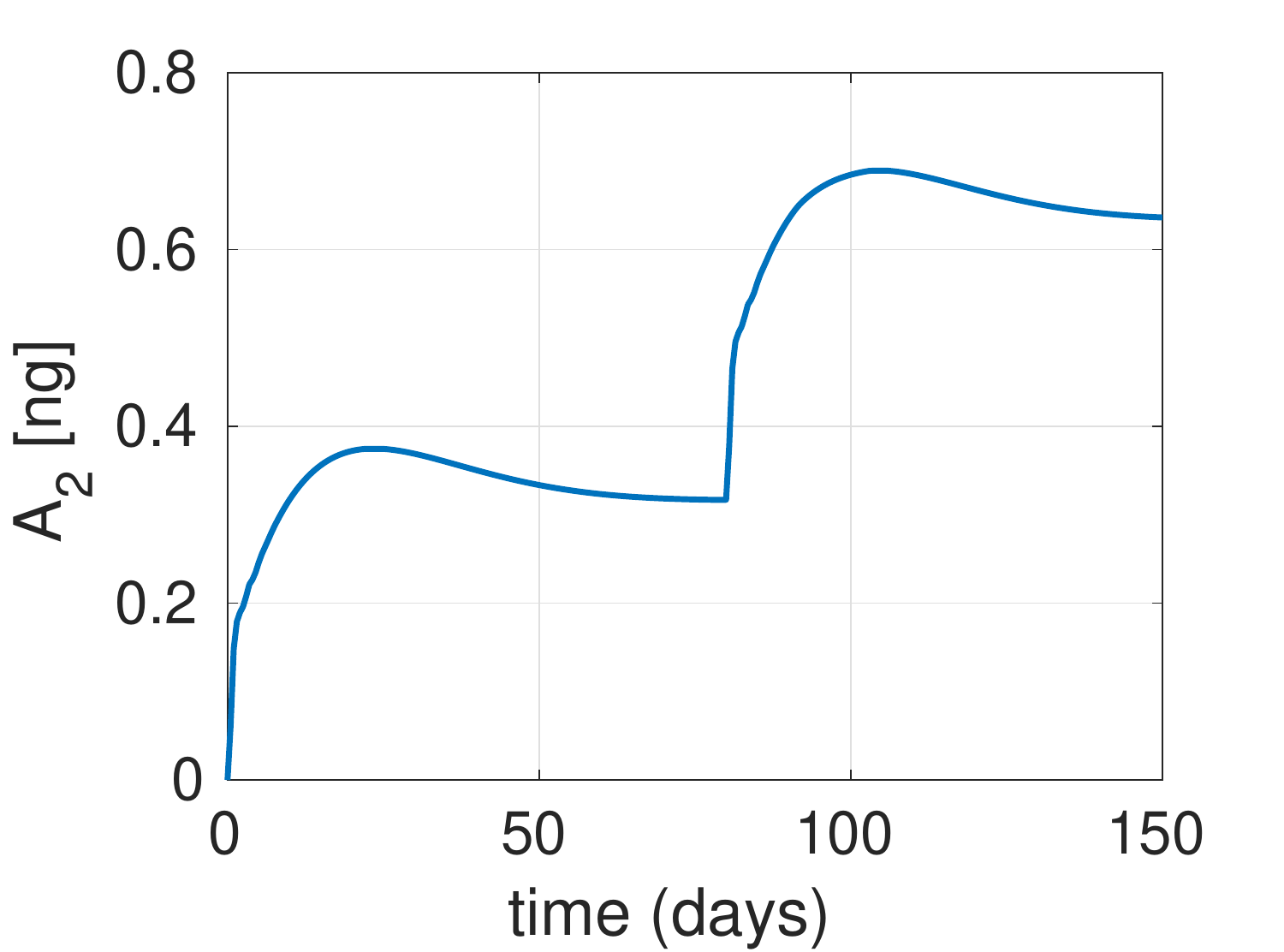}{}%
\includegraphics[width=0.325\columnwidth]{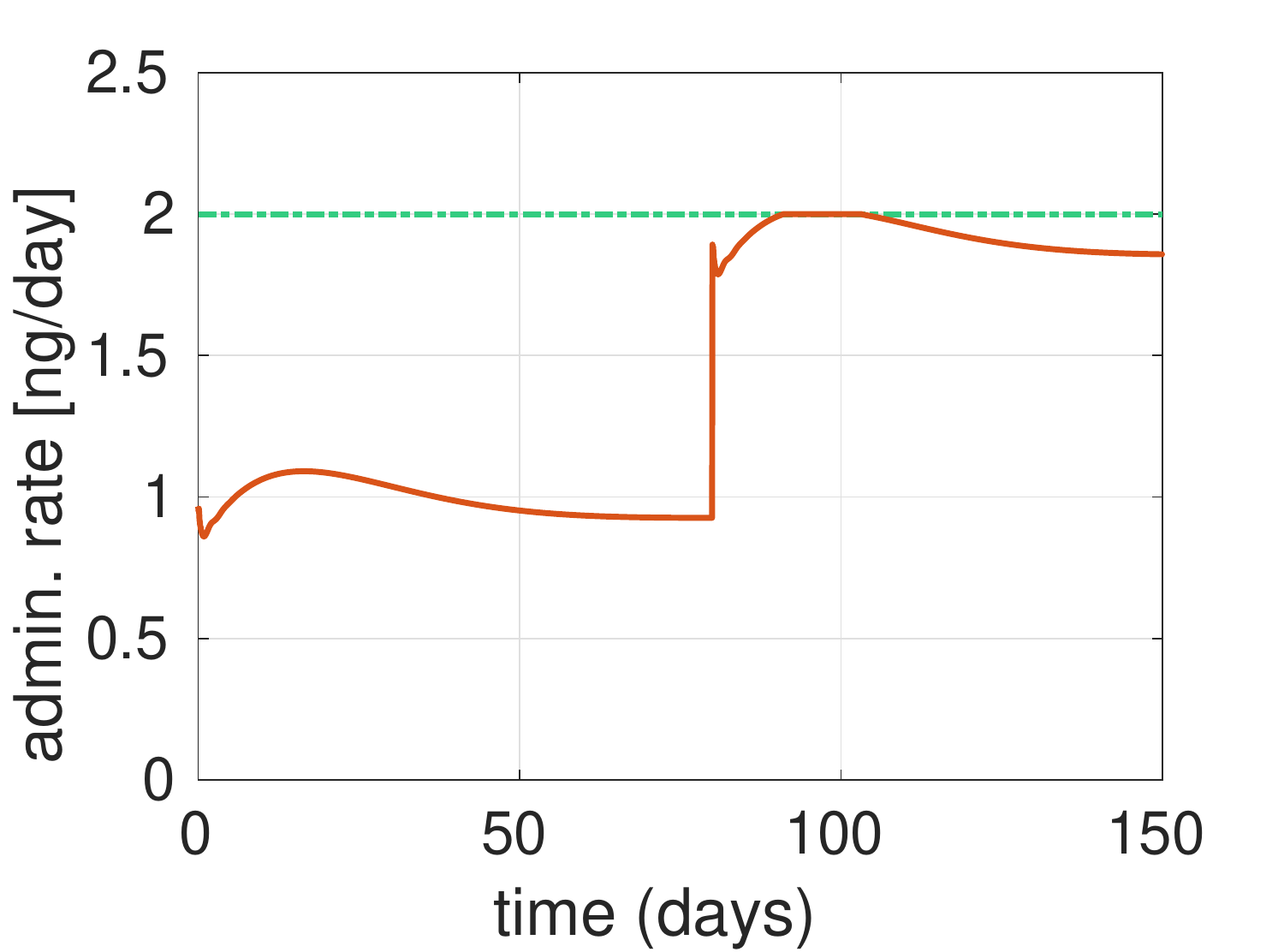}
  \caption{ (Left) Controlled amount of drug in the central compartment ($A_1$) along with the prescribed reference signal and the upper bound on $A_1$, (Middle) amount of drug in the peripheral compartment, (Right) Administration rate determined by the offset-free MPC controller.}
    \label{fig:fracResponse}
\end{figure}

The quantity of amiodarone in the central compartment, $A_1$,  is shown in Figure~\ref{fig:fracResponse} (left) where it can be seen that it tracks the prescribed set-point while it does not exceed the constraint of $\unit[1.03]{ng}$.
The input produced by the controller as the solution of the optimization problem $\mathsf{P}_{\mathrm{mpc}}$ is shown in Figure~\ref{fig:fracResponse} (right) and it may be observed that the administration rate does not exceed the maximum allowed value of $\unitfrac[2.0]{ng}{day}$.
In Figure~\ref{fig:fracResponse} (middle) we show the amount of drug in the peripheral/tissues compartment.

\subsection{Sensitivity analysis}
In order to assess the sensitivity of the closed-loop system to the 
exact knowledge of the pharmacokinetic parameters $k_{12}$, $k_{21}$, 
$k_{10}$ and $\alpha$ we performed closed-loop simulations against perturbed models using the MPC designed for the exact model. In particular 
we assumed a perturbation of $+10\%$ or $-10\%$ on each one of the pharmacokinetic parameters. 

In each closed-loop simulation, we assumed the same scenario as in the nominal case, i.e., the treating 
physician selects a set-point equal to $\unit[0.5]{ng}$ for the first $80$ days 
and increases it to $\unit[1.0]{ng}$ for the next $70$ days. 
We also assumed the same upper bounds on the administration rate 
and the quantity of amiodarone in the central compartment, $A_1$. In order to 
evaluate the performance of the drug administration course, we introduced the 
index
\begin{align}
 J = \frac{1}{N_u}\sum_{k=0}^{N_u-1}(y_k - r_k)^2 + u_k^2,
\end{align}
where $N_u$ is the number of simulation points (here $N_u=3000$).
This performance indicator quantifies and encodes the average tracking error and the 
average administration rate, both of which are desirable to be low.

\begin{figure}[ht!]
\centering
\includegraphics[width=0.325\columnwidth]{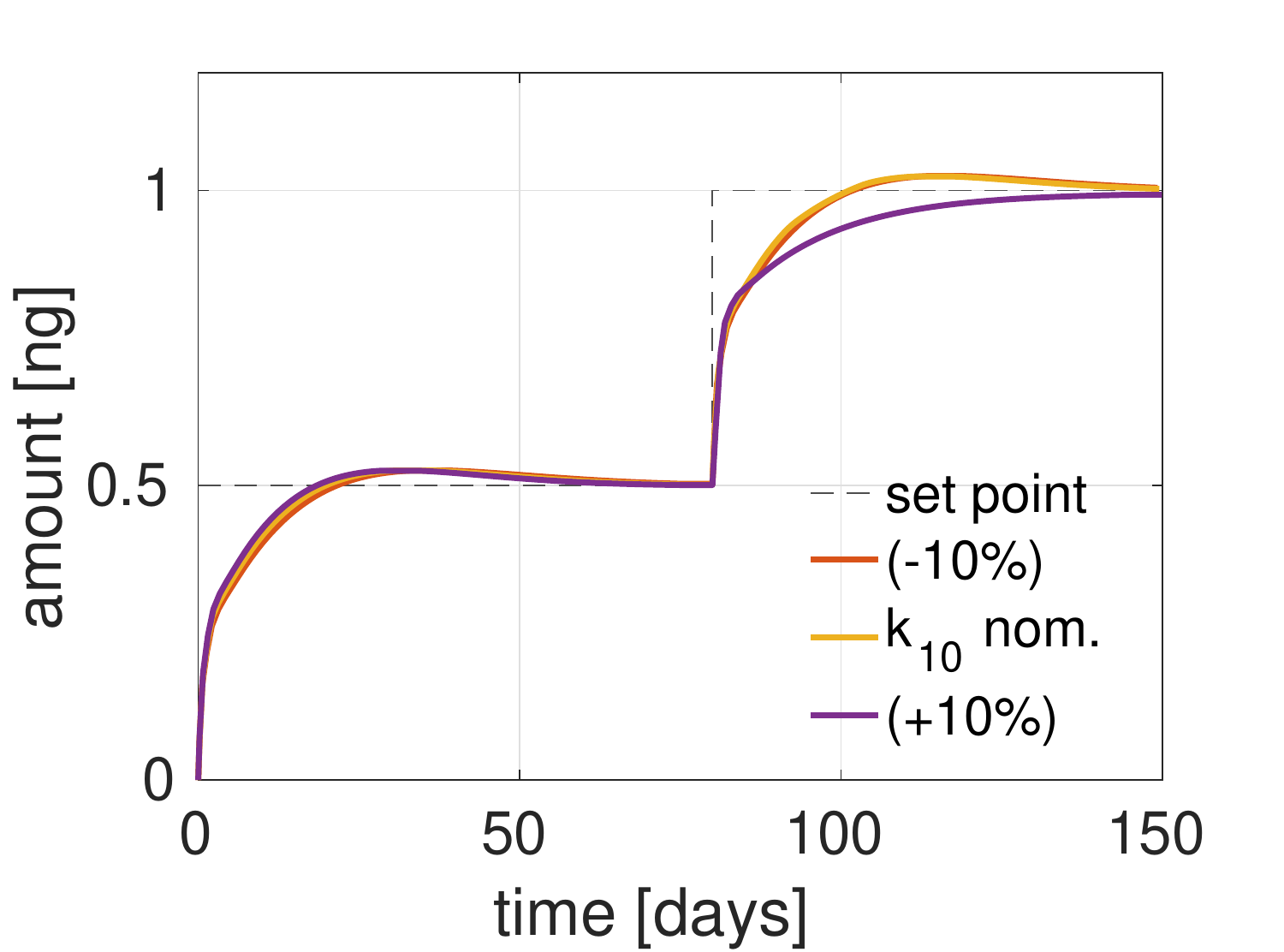}{}%
\includegraphics[width=0.325\columnwidth]{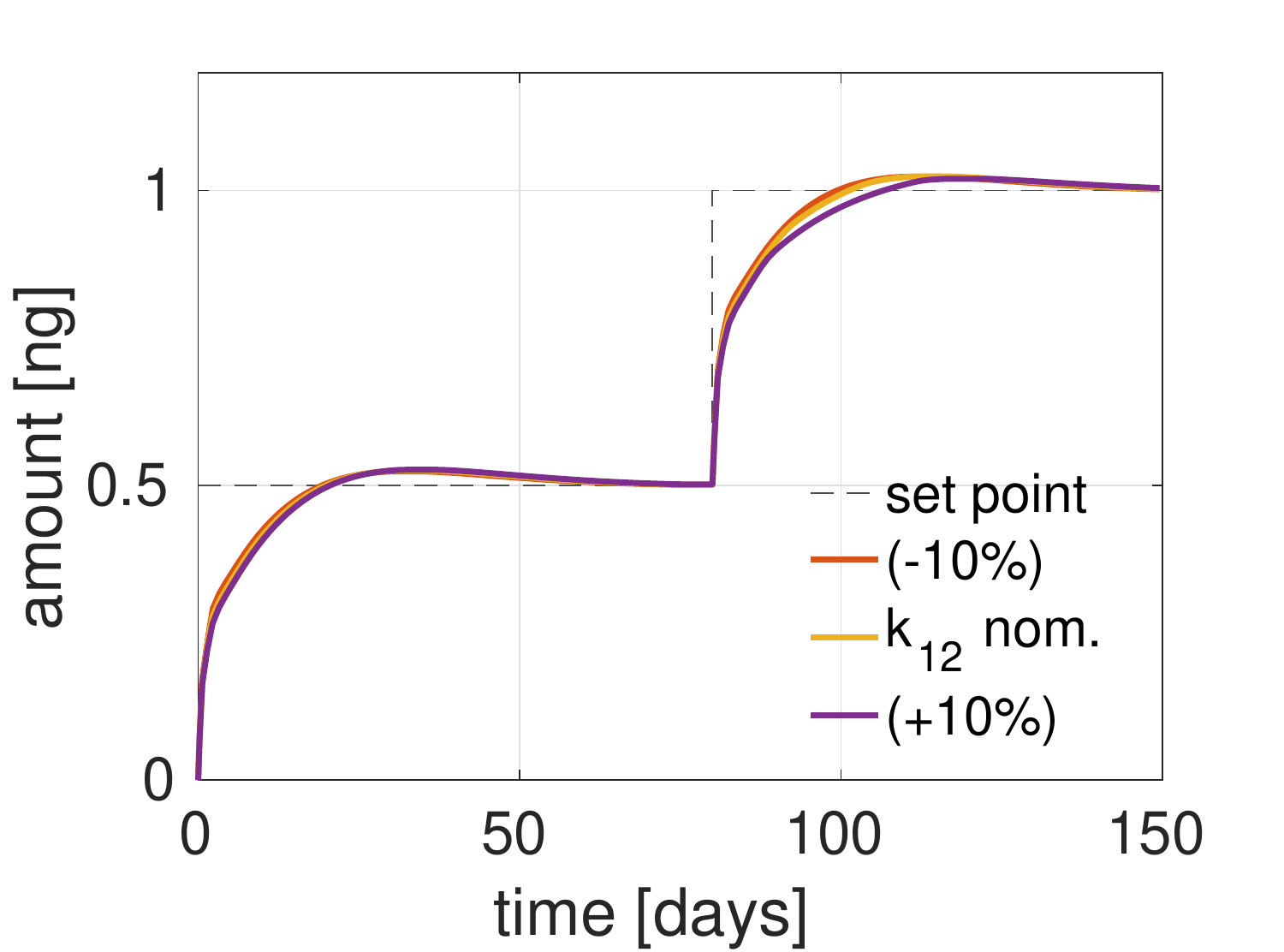}{}%
\includegraphics[width=0.325\columnwidth]{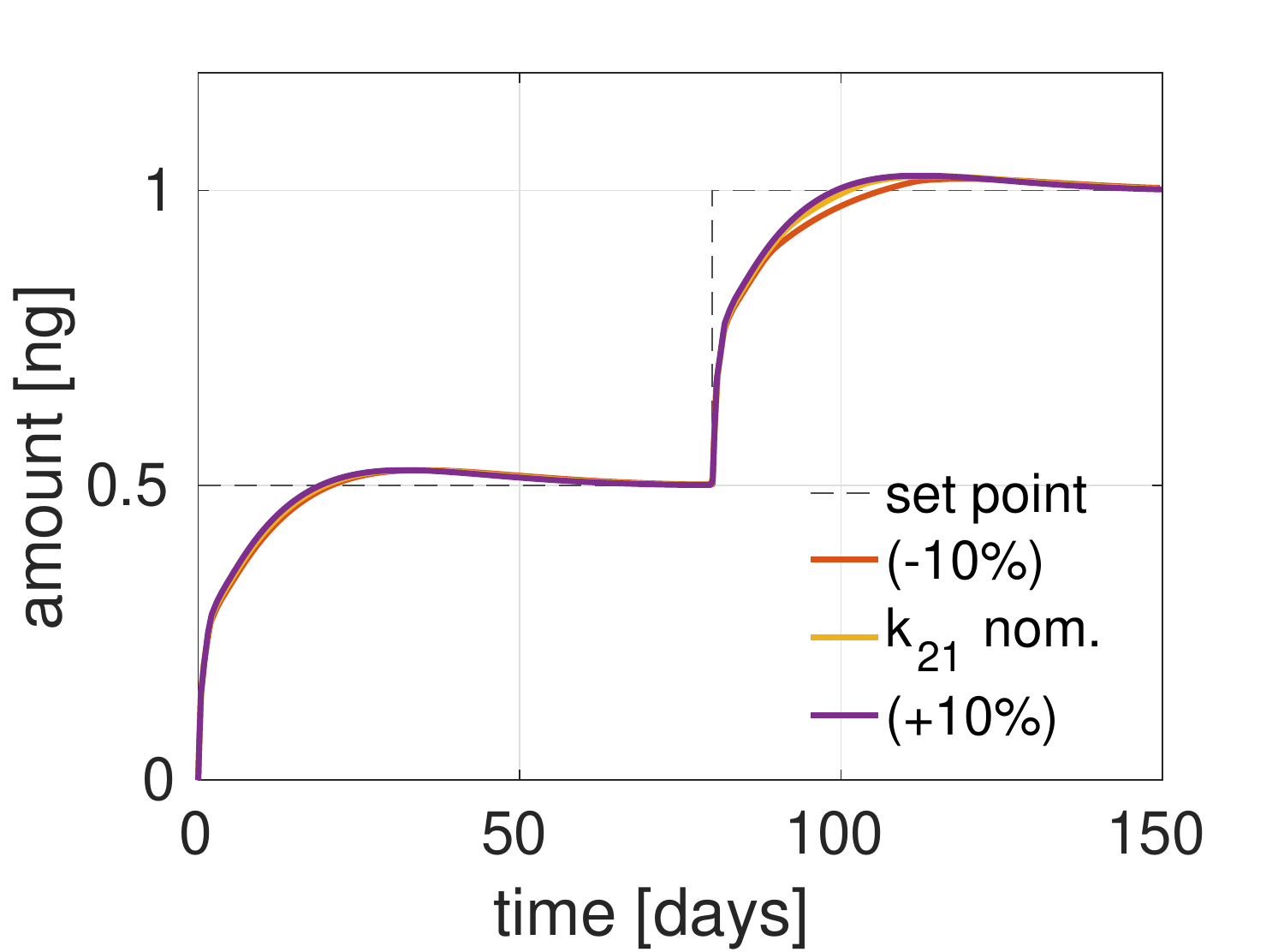}{}\\
\includegraphics[width=0.325\columnwidth]{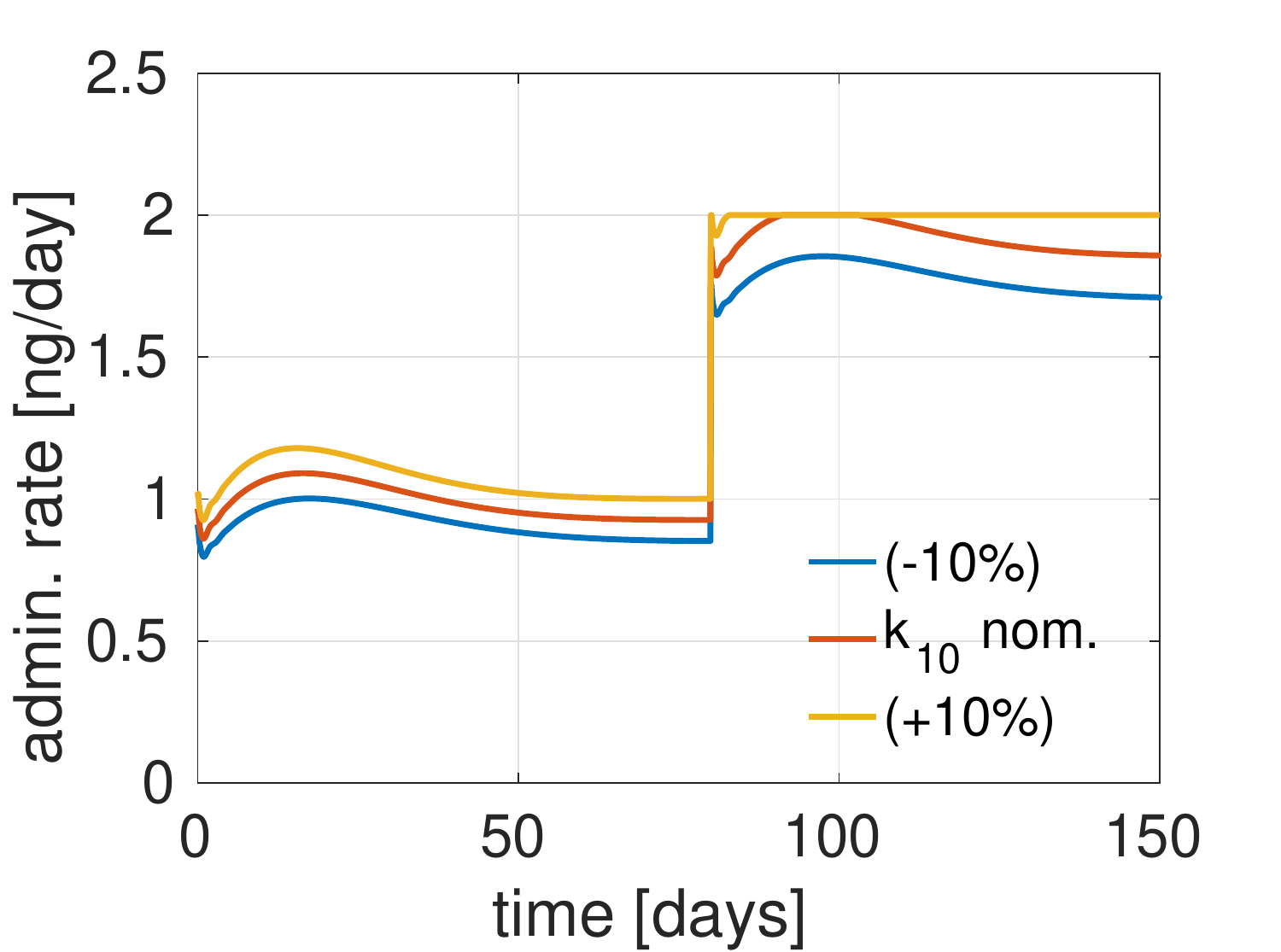}{}%
\includegraphics[width=0.325\columnwidth]{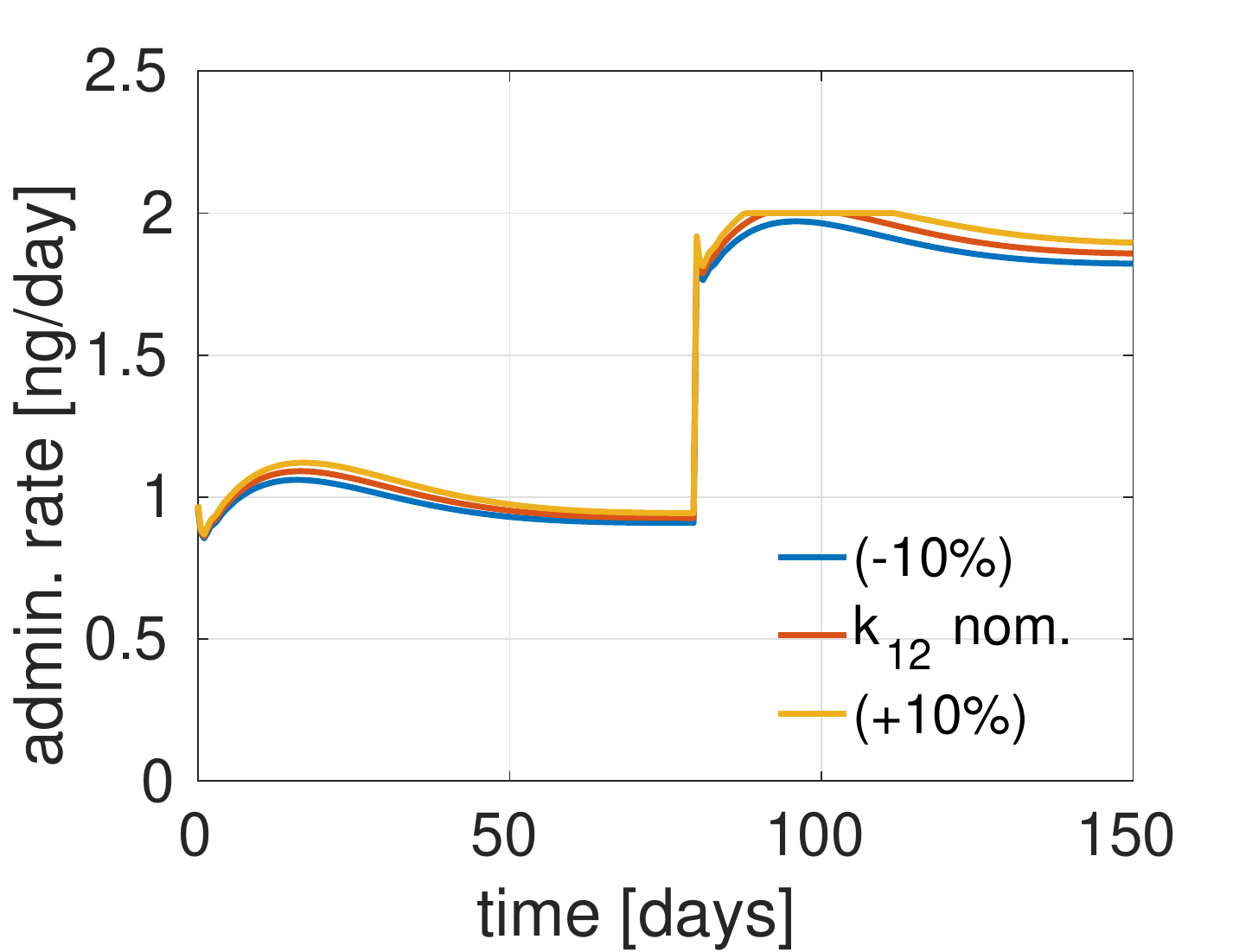}{}%
\includegraphics[width=0.325\columnwidth]{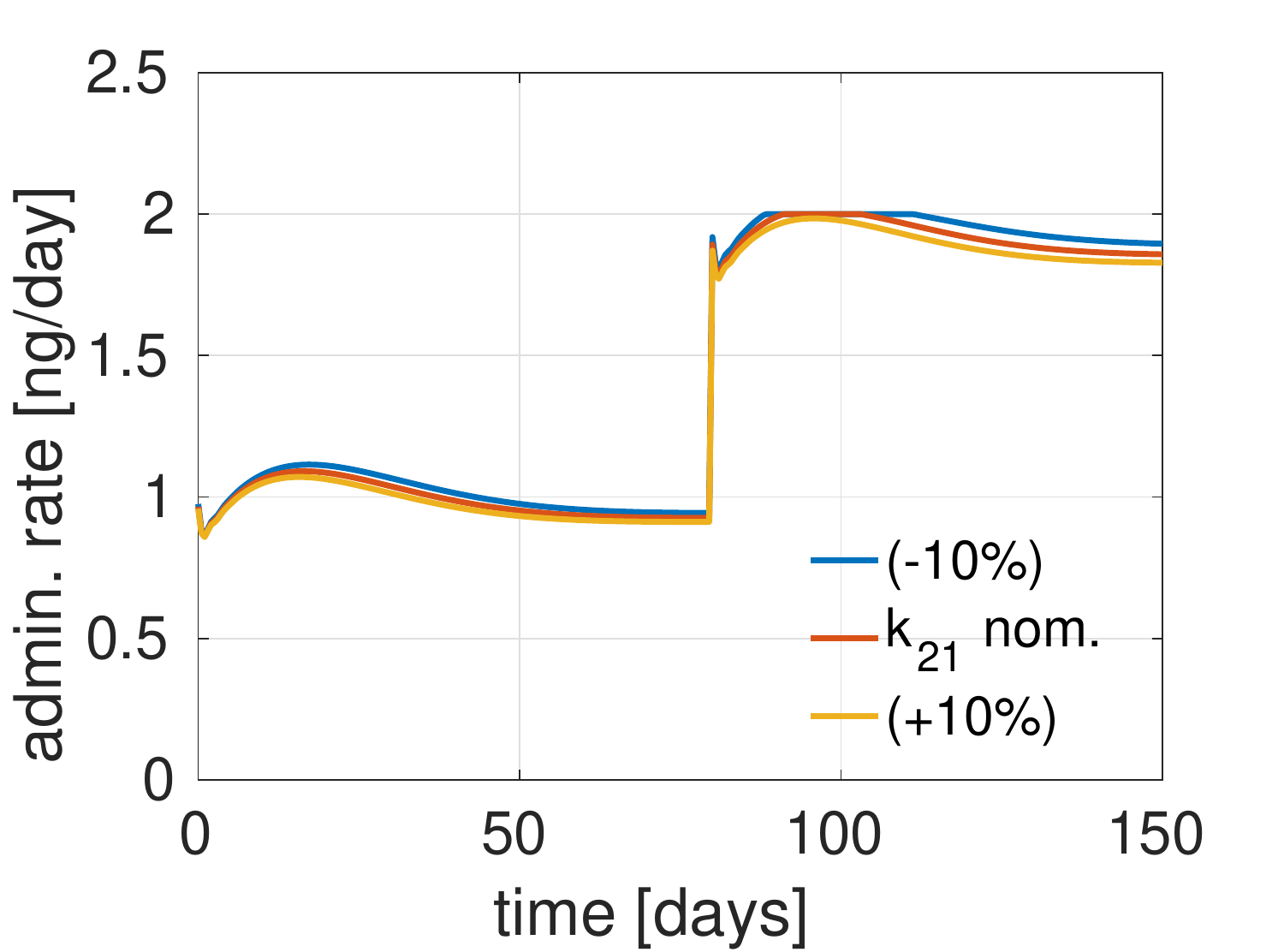}{}\\
\includegraphics[width=0.325\columnwidth]{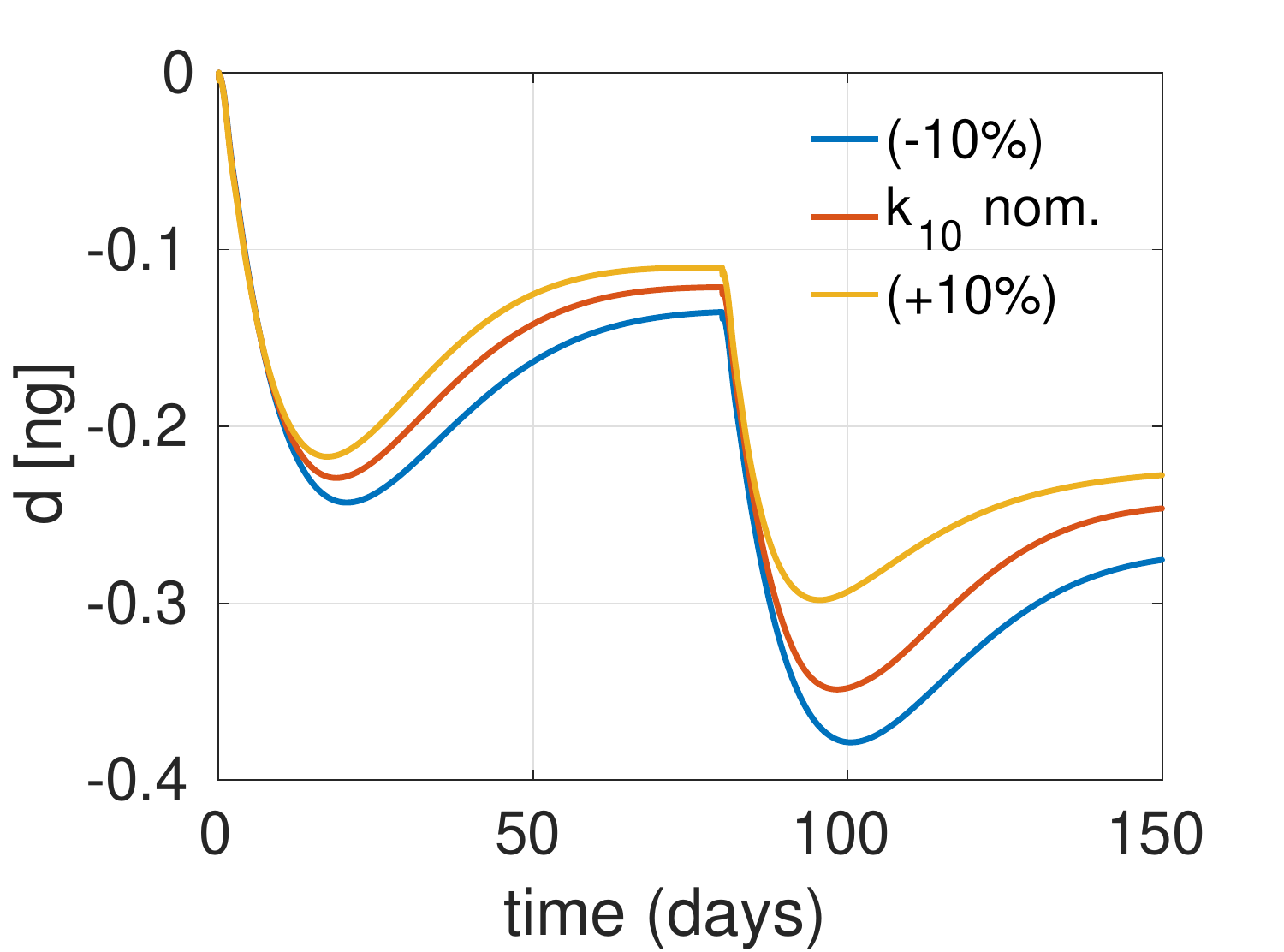}{}%
\includegraphics[width=0.325\columnwidth]{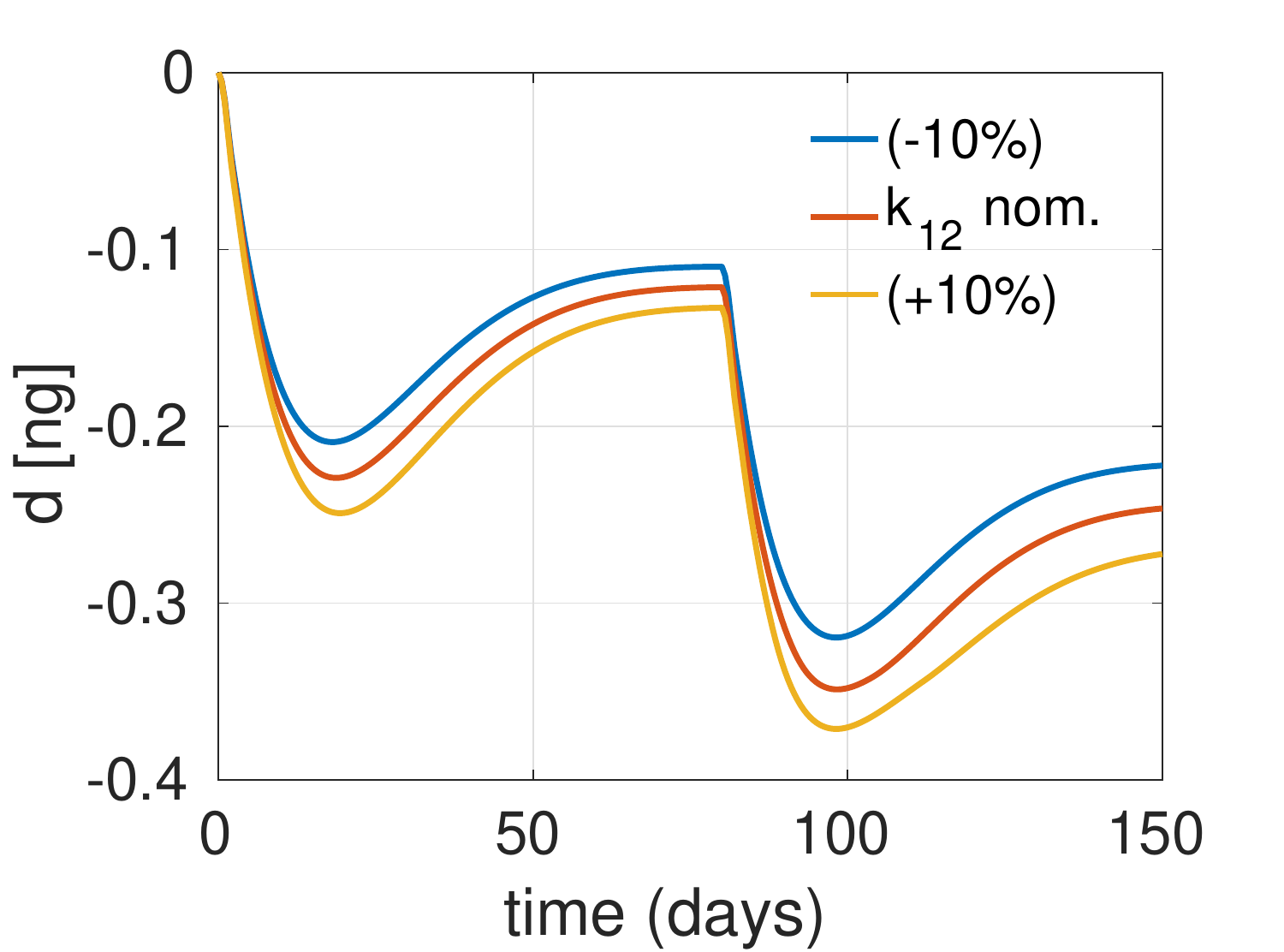}{}%
\includegraphics[width=0.325\columnwidth]{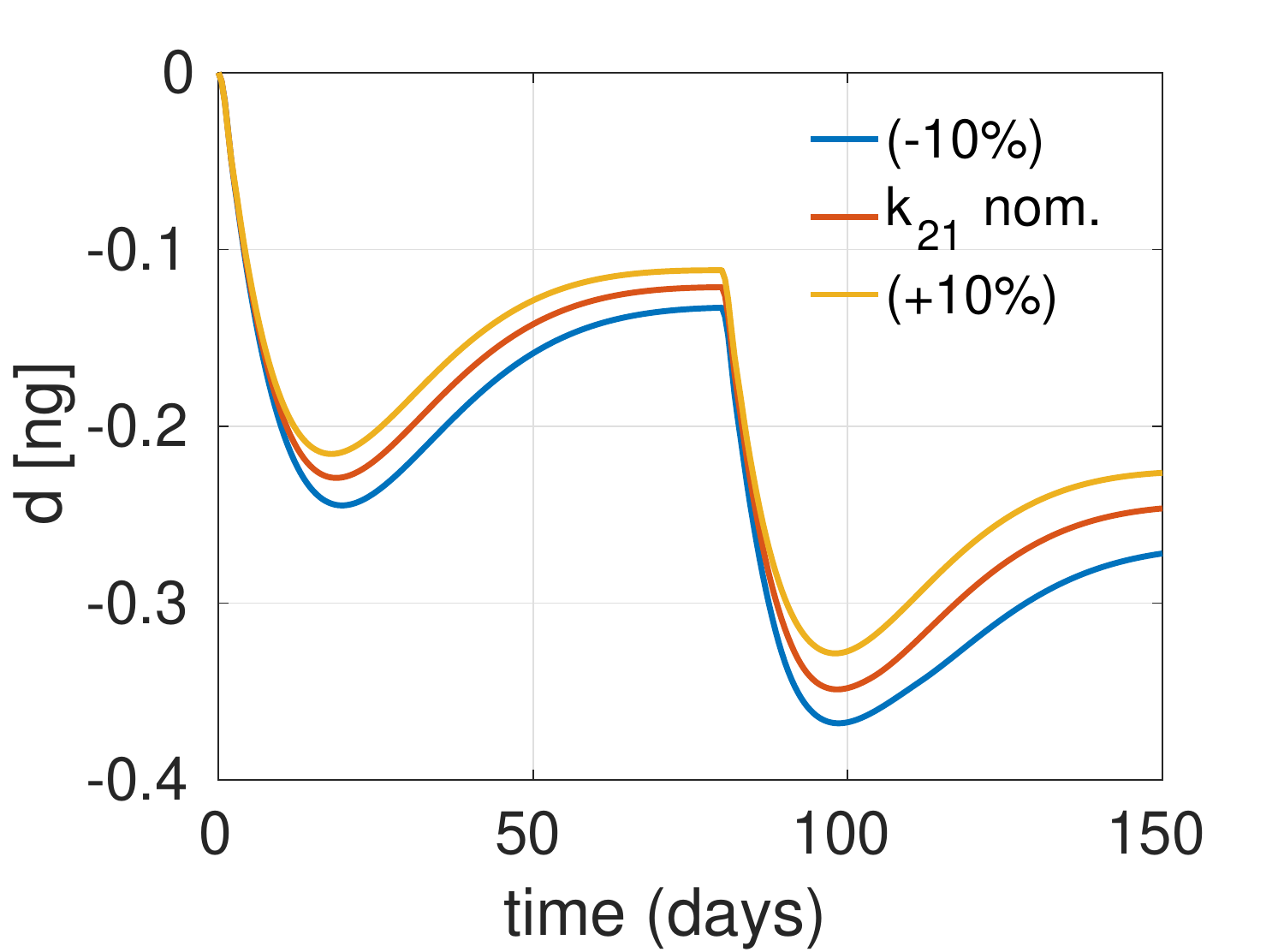}{}
  \caption{(First row) Controlled response of $A_1$, the amount of drug in the central compartment following two changes in the set-point and sensitivity to $k_{10}$ (left), $k_{12}$ (middle) and $k_{21}$ (right). (Second row) Drug administration rate $u_k$ computed by the offset-free MPC controller with $N=60$ and $\nu=25$. (Third row) Corresponding disturbance estimates, $\hat{d}_k$.}
    \label{fig:sensitivity-PK}
\end{figure}

\begin{figure}[ht!]
\centering
\includegraphics[width=0.32\columnwidth]{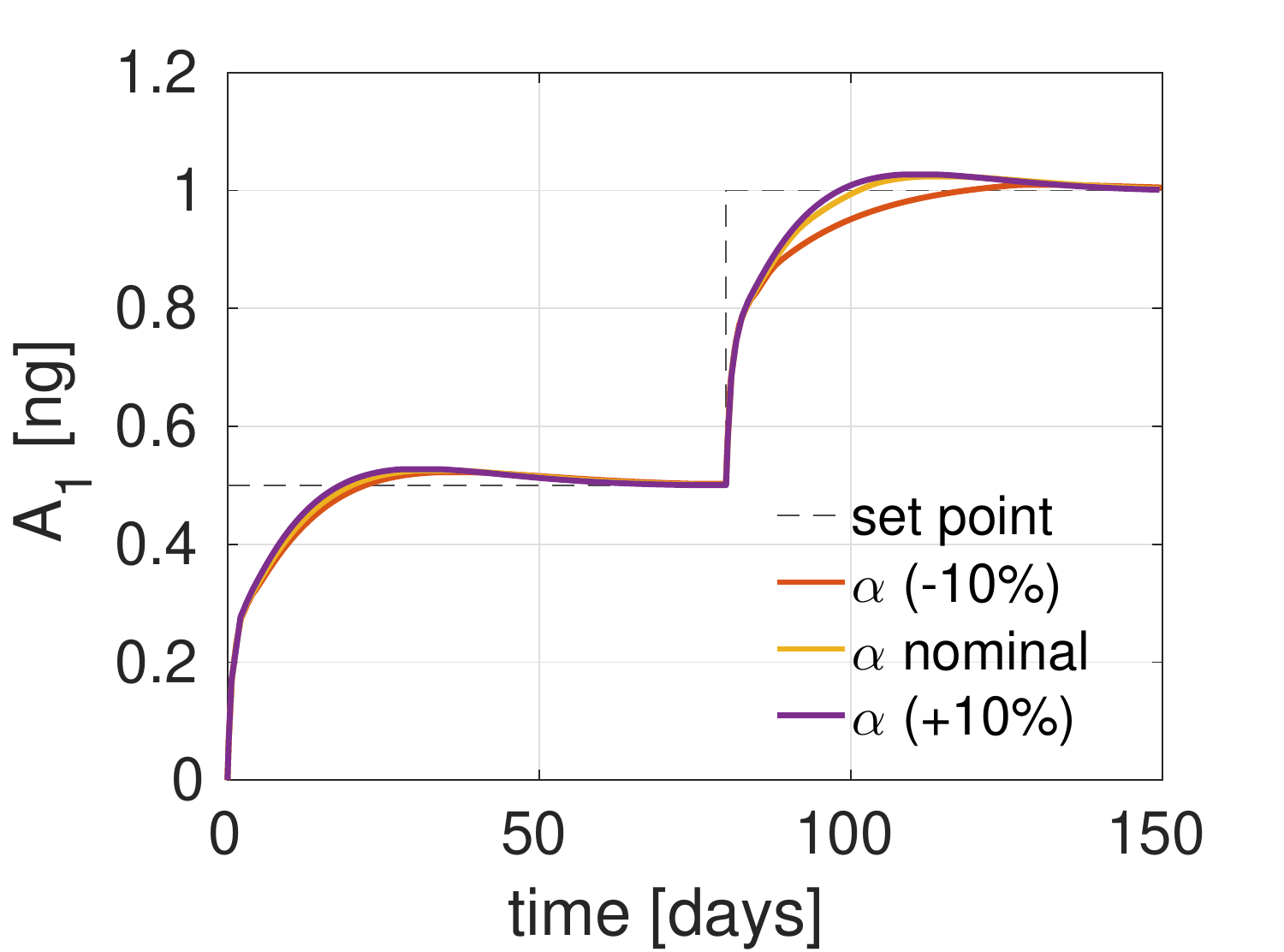}{}%
\includegraphics[width=0.32\columnwidth]{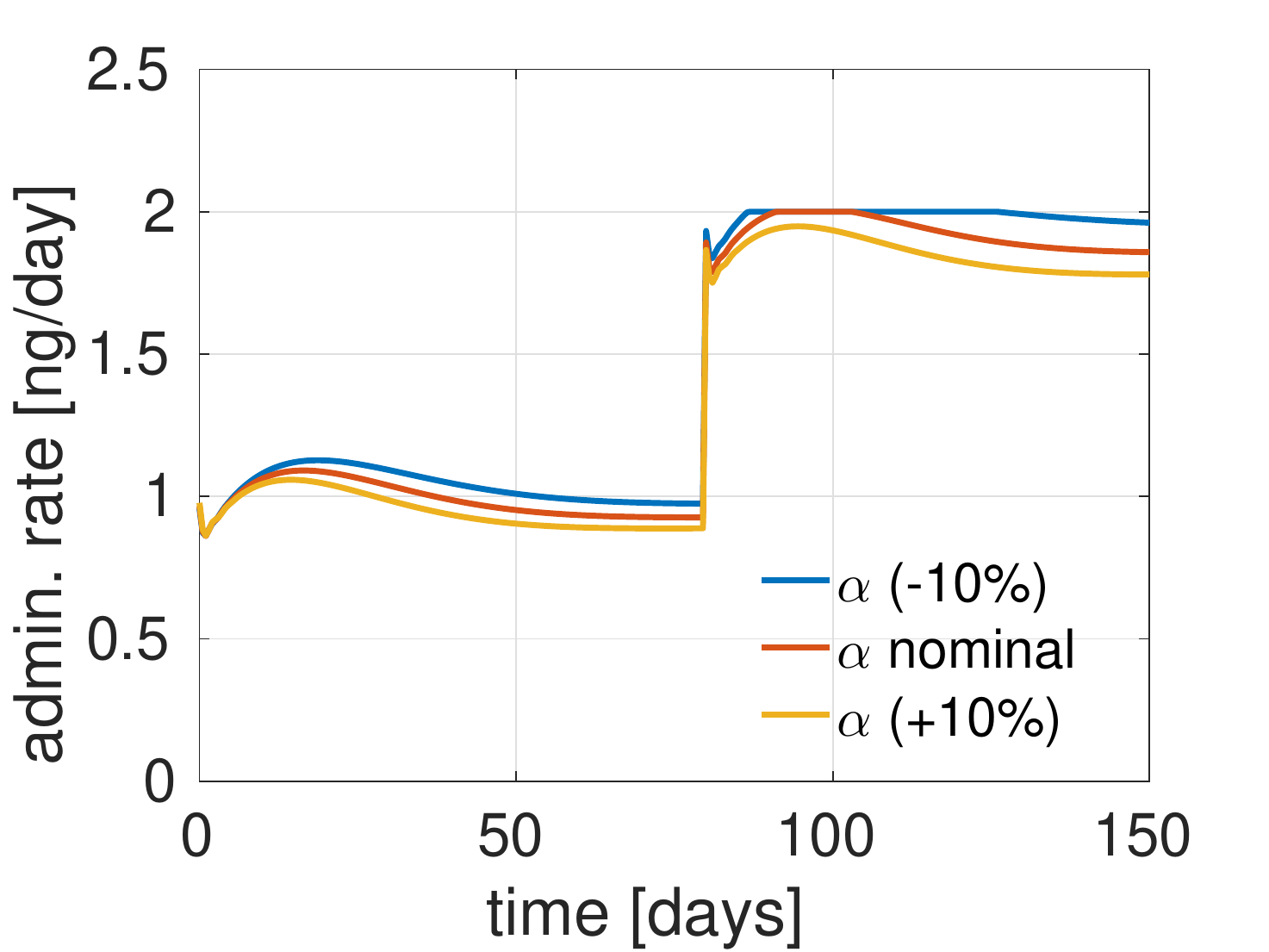}{}%
\includegraphics[width=0.32\columnwidth]{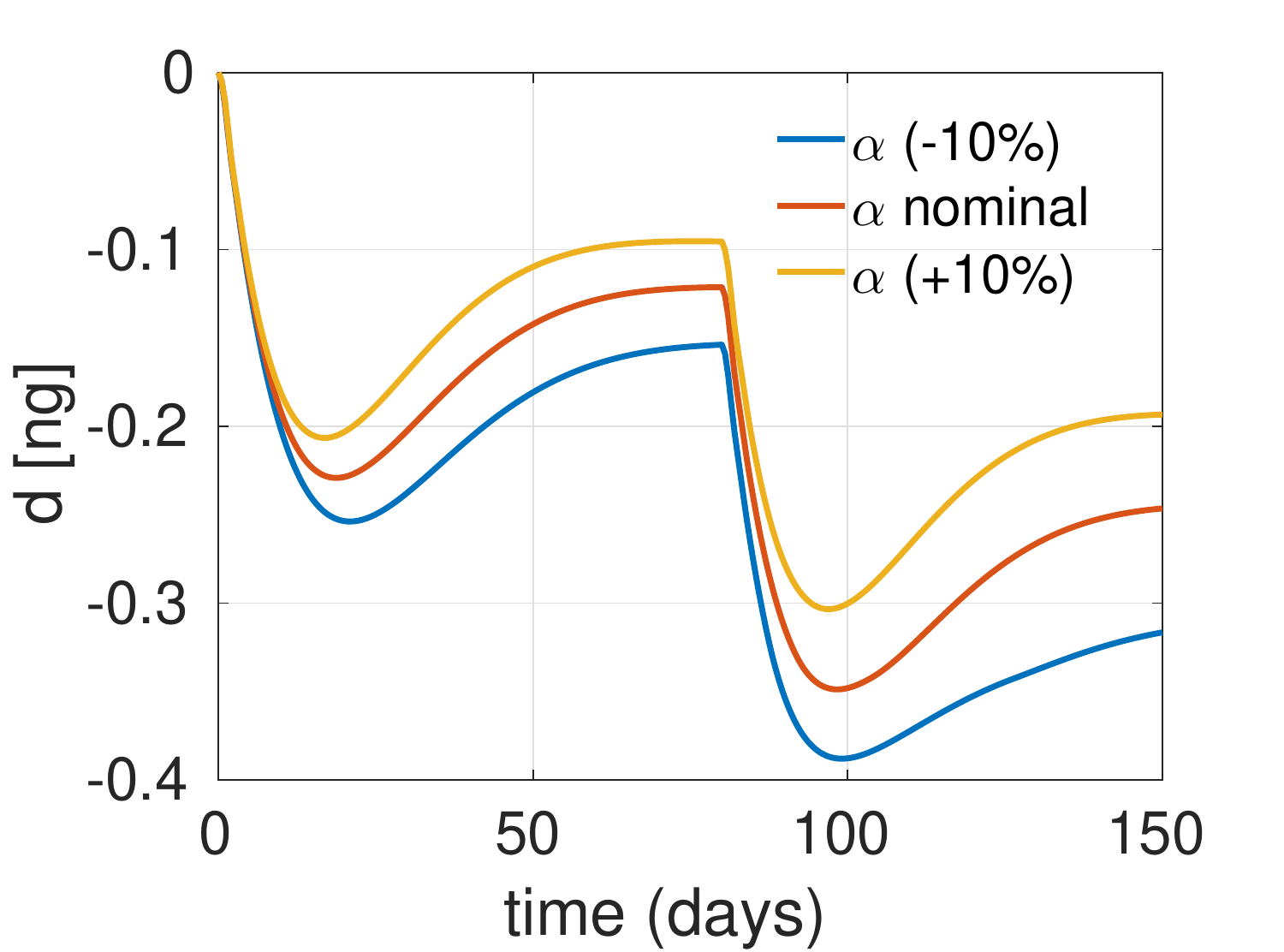}
  \caption{ (Left) Controlled response of $A_1$ and (Middle) administration rate computed by the offset-free MPC controller with $N=60$ and $\nu=25$ for different values of $\alpha$. (Right) Estimated disturbance, $\hat{d}_k$, for different values of $\alpha$.}
    \label{fig:sensitivity-alpha}
\end{figure}

  Figures ~\ref{fig:sensitivity-PK} and  ~\ref{fig:sensitivity-alpha} present the input-output responses and the disturbance estimates~$\hat{d}_k$ and Table ~\ref{tab:J1}  shows the value of the performance index for the nominal case and all eight perturbations. It seems that the closed-loop system is more sensitive 
  to errors in the clearance rate constant $k_{10}$, however its effect 
  is rather limited. Evidently, in all simulation results, the impact of parametric uncertainty
  is not significant and 
  the offset-free control scheme produces successful control actions in tracking the prescribed set-point changes. Overall, despite the fact that the exact values 
  of $k_{10}$, $k_{12}$ and $k_{21}$ are not available to the controller, the prescribed 
set-point on $A_1$ is reached and the constraints on $A_1$ and $u$ are satisfied at 
  all time instants.

\begin{figure}[ht!]
\centering
\includegraphics[width=0.32\columnwidth]{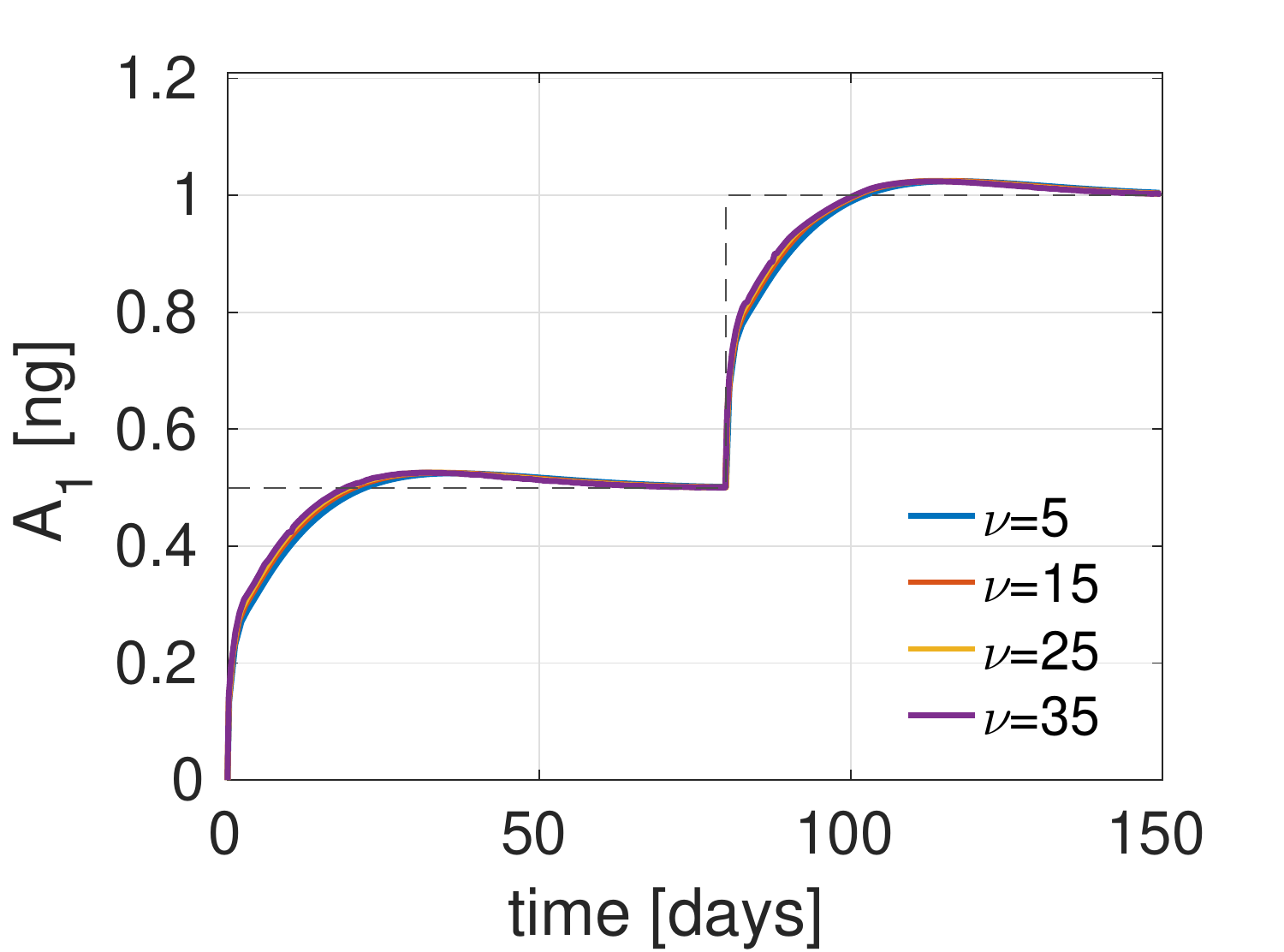}{}%
\includegraphics[width=0.32\columnwidth]{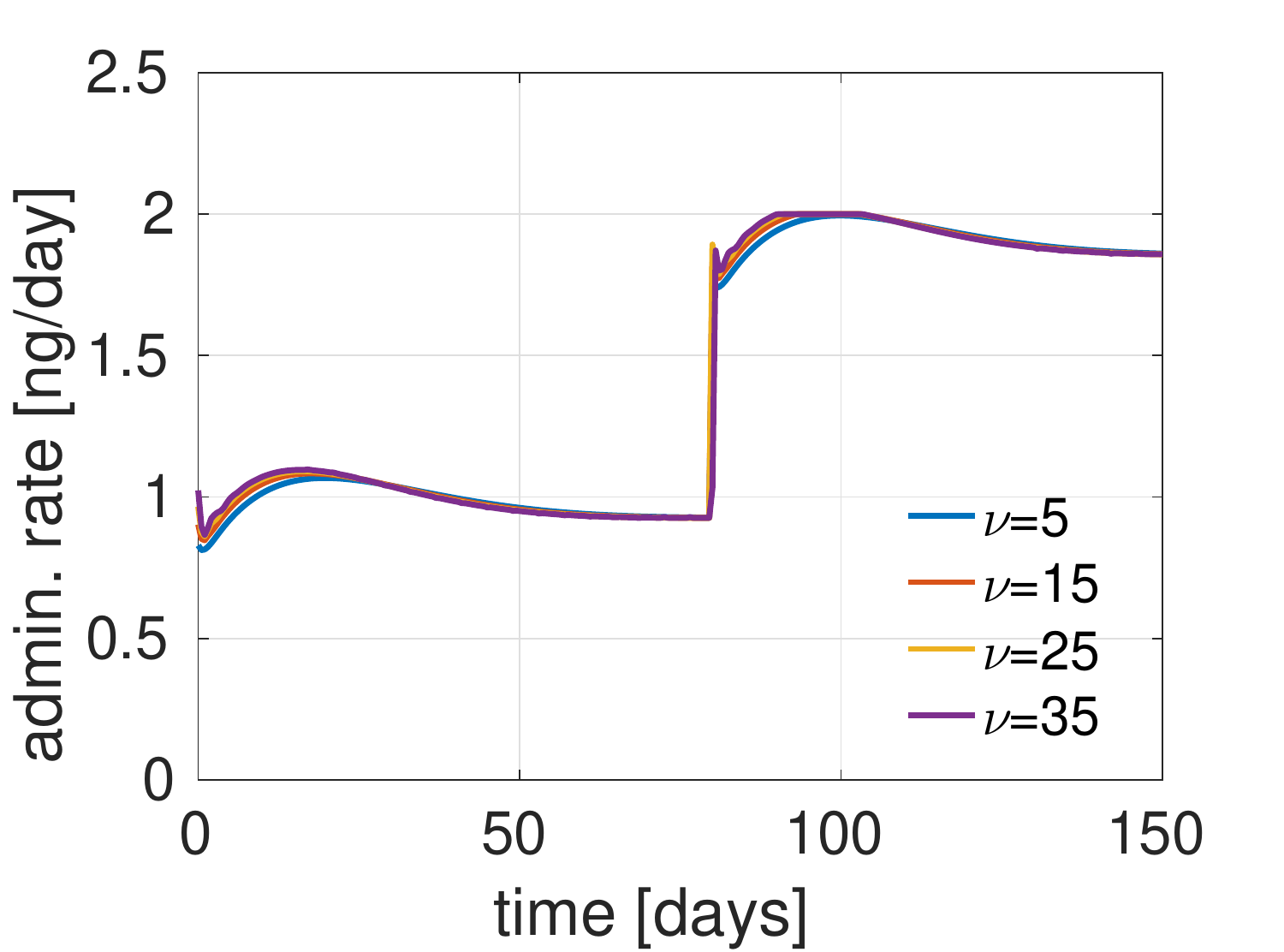}{}%
\includegraphics[width=0.32\columnwidth]{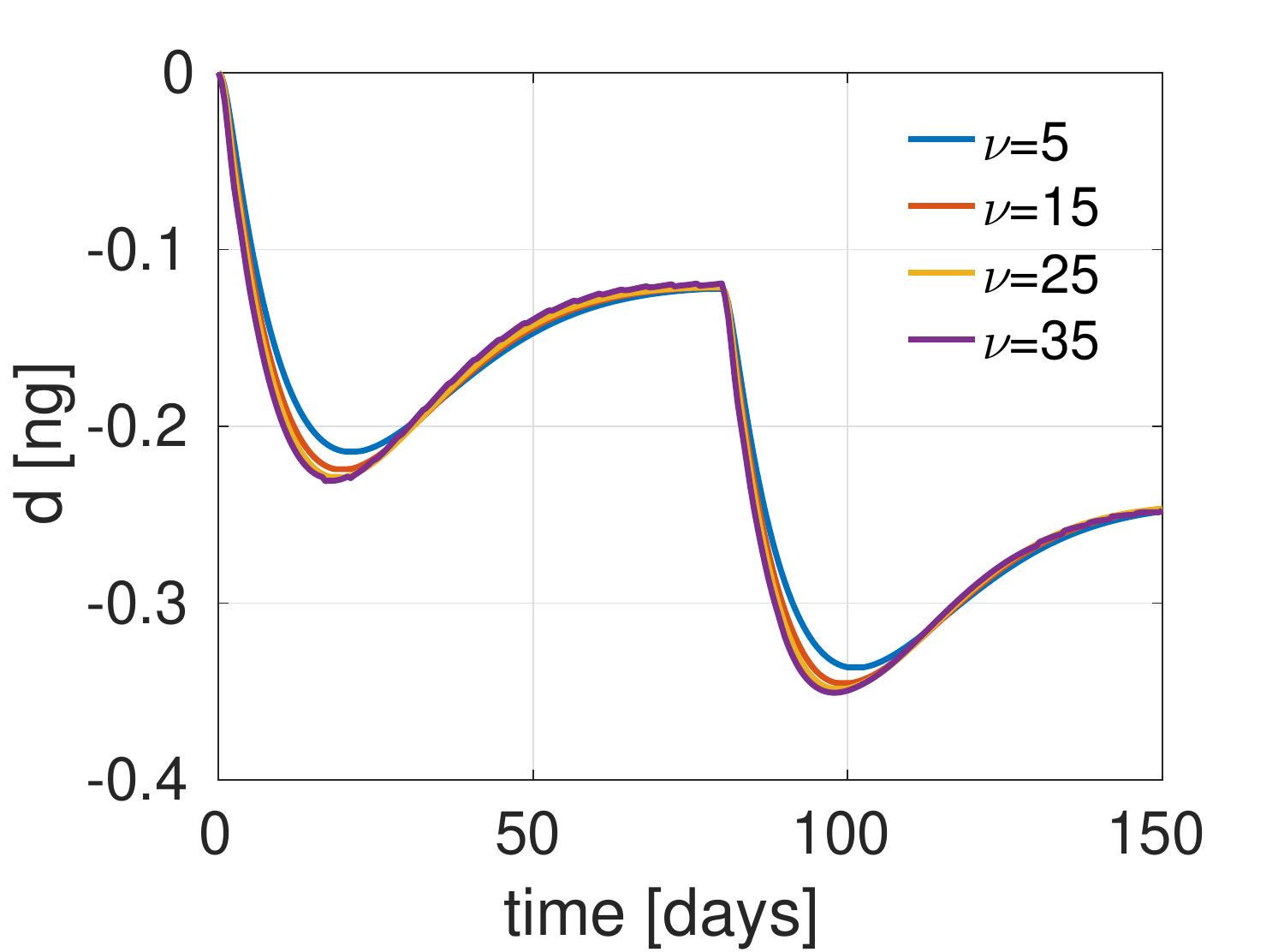}
  \caption{Controlled response of $A_1$ (left) and administration rate (middle) computed by the offset-free MPC controller with $N=60$ and for different memory lengths $\nu$. (Right) Disturbance estimates, $\hat{d}_k$, for different values of $\nu$.}
    \label{fig:sensitivity-nu}
\end{figure}

  In the following simulations, we investigated the effect of the value of the memory length $\nu$ on the performance of the system. Again, we ran the same scenario with the nominal case concerning the set point changes and the constraints of the input and state parameters using memory lengths equal to $5$, $15$, $25$ and $35$. As we discussed in 
  Section~\ref{sec:fin-dim-approx}, larger values of $\nu$ lead to more 
  accurate models. As one may observe in Figure~\ref{fig:sensitivity-nu}
  and Table~\ref{tab:J2}, the performance is only slightly improved by increasing the value of a memory length $\nu$ and in general the closed-loop performance is not sensitive to the selection of $\nu$. For  $\nu=25$ the response curves
  are almost indistinguishable from those with $\nu=35$.

\begin{table}[ht!]
\centering
\begin{tabular}{ l | c | c | c }
  \hline			
   & $-10\%$ & Nominal & $+10\%$\\ \hline
  $k_{10}$ & 9.9628 & 11.6345 & 12.7328 \\
  $k_{12}$ & 11.1457 &  & 12.0466 \\
  $k_{21}$ & 12.0169 & & 11.2554 \\
  $\alpha$ & 12.4379 & & 10.7575 \\
  \hline  
\end{tabular}
\caption{Values of the performance indicator $J$ for perturbed values of pharmacokinetic parameters for $N=60$ and $\nu=25$.}
\label{tab:J1}
\end{table}
\begin{table}[ht!]
\centering
\begin{tabular}{ l | c }
  \hline			
  $\nu$ & $J$\\ \hline  
  5 &   11.4950 \\
  15 &  11.5998 \\
  25 &  11.6345 \\
  35 &  11.5931 \\
  \hline  
\end{tabular} 
\caption{Values of the performance indicator $J$ for different values of the memory length for $N=60$.}
\label{tab:J2}
\end{table}

Overall, we see that parametric errors can be well attenuated and have little effect on 
the performance of the closed-loop controlled system. As one may observe in Figures~\ref{fig:sensitivity-PK},
\ref{fig:sensitivity-alpha} and~\ref{fig:sensitivity-nu}, the disturbance estimates~$\hat{d}_k$ reconstruct 
the model-system mismatch and tracking offset is avoided.

 
\subsection{Computational aspects}
All simulations were run in {\textsc MATLAB}\textsuperscript{\textregistered} (version 2014b) on a Macbook Pro (version 10.12.4, \unit[2.66]{GHz} Intel Core 2 Duo, \unit[4]{GB} RAM) running Mac OS Sierra. The MPC problem was formulated as a constrained quadratic program which can be solved very efficiently in real time. The solver used for solving the optimization problem was {\textsc{Mosek}}~\cite{mosek}. The average runtime to solve the MPC problem was \unit[1.3]{s} and the maximum runtime was \unit[2.3]{s} for a memory of $\nu=35$. 

\section{Conclusions}
In this paper an offset-free tracking MPC scheme was proposed for 
constrained fractional-order systems and it was demonstrated that it 
is resilient to modeling errors. 

The proposed tracking methodology is based on a discrete-time 
finite-dimension approximation of the infinite-dimensional fractional dynamics 
and can handle constraints on both state and input variables. 
The approximation used in this work hinges on the 
Gr\"unwald-Letnikov fractional-order derivative, it is a time-domain 
approximation and, unlike frequency-domain approximations, leads to 
a bounded approximation error with known error bounds.
Combined with the design of an observer, the proposed offset-free MPC approach 
is compatible with the availability of partial state information, 
which makes it suitable for real control problems. The augmented state observer which is used in this work produces
estimates both for the system state and a disturbance term which 
encompasses modeling errors. Asymptotically, the observer reconstructs 
both the states and the disturbance term leading to offset-free 
tracking despite possible modeling errors.


The merits of the proposed MPC scheme were demonstrated on the problem 
of intravenous administration of amiodarone; a drug which follows fractional-order 
pharmacokinetics. Using blood stream measurements only, the closed-loop
system, which combines an augmented state observer with an MPC controller,
leads to zero-offset tracking and is resilient to modeling errors. 
This characteristic of the proposed methodology deems it suitable for 
drug administration applications where the actual pharmacokinetic parameters 
of the patient are unknown.
We demonstrated that the sensitivity of the closed-loop 
system to the errors on the pharmacokinetic parameters is rather low.
Additionally, a memory length of $\nu=25$ is adequate for controller 
design purposes leading to a computationally tractable formulation.
The resulting optimization problems are solved very efficiently
and do not pose any limitations to the applicability of the controller.

The theory and the results presented in this paper demonstrate that the proposed MPC scheme inherits the key advantages of offset-free MPC to fractional-order systems in terms of satisfying the system constraints and achieving zero steady-state error even in presence of modeling uncertainties. Future research directions will examine extensions of this scheme to formulate robust or stochastic MPC strategies and take into account time delays.




\section*{References}
\bibliographystyle{model1-num-names}
\bibliography{sample.bib}







\end{document}